\documentclass[leqno,11pt]{amsart}
\usepackage{amsmath,amscd}
\usepackage{epsfig,graphics}
\usepackage{amssymb,latexsym}
\usepackage{mathrsfs}
\usepackage{eucal}

\newtheorem{thm}{\bf Theorem}[section]

\newtheorem{prop}[thm]{\bf Proposition}
\newtheorem{cor}[thm]{\bf Corollary}
\newtheorem{lem}[thm]{\bf Lemma}
\newtheorem{rem}[thm]{\bf Remark}

\newcommand{\A}{\mathcal{A}}
\newcommand{\cA}{\mathscr{A}}
\newcommand{\cD}{\mathscr{D}}
\newcommand{\B}{\mathbf{B}}
\newcommand{\cB}{\mathbf{B}}

\newcommand{\cP}{\mathscr{P}}
\newcommand{\R}{\mathscr{R}}
\newcommand{\pf}{\noindent{\bfseries Proof. }}
\newcommand{\ov}{\overline}

\newcommand{\E}{{\mathcal{E}}}

\newcommand{\bM}{{\mathbf{M}}}
\newcommand{\F}{\mathcal{F}}

\newcommand{\gl}{\mathfrak{gl}}
\newcommand{\Z}{\mathbb{Z}}
\newcommand{\C}{\mathbb{C}}
\newcommand{\h}{\mathfrak{h}}
\newcommand{\te}{\widetilde{e}}
\newcommand{\tf}{\widetilde{f}}
\newcommand{\g}{\mathfrak{g}}
\newcommand{\td}{\widetilde}
\newcommand{\tE}{\widetilde{E}}
\newcommand{\tF}{\widetilde{F}}
\newcommand{\K}{\mathscr{K}}

\numberwithin{equation}{section}

\begin{document}
\title[Littlewood-Richardson rule of extremal weight crystals]
{Crystal duality and Littlewood-Richardson rule of extremal weight
crystals}
\author{JAE-HOON KWON}
\address{Department of Mathematics \\ University of Seoul   \\  Seoul 130-743, Korea }
\email{jhkwon@uos.ac.kr }

\thanks{This work was supported by LG Yonam foundation.}

\begin{abstract}
We consider a category of $\gl_\infty$-crystals, whose objects are
disjoint unions of extremal weight crystals of non-negative
level with certain finite conditions on the multiplicity of connected components. We show
that it is a monoidal category under tensor product of crystals and
the associated Grothendieck ring is anti-isomorphic to an Ore
extension of the character ring of integrable lowest weight
$\gl_\infty$-modules with respect to derivations shifting the
characters of fundamental modules. A Littlewood-Richardson rule of
extremal weight crystals with non-negative level is described
explicitly in terms of classical Littlewood-Richardson coefficients.
\end{abstract}

\maketitle

\section{Introduction}
Let $U_q(\frak{g})$ be the quantized enveloping algebra associated
with a symmetrizable Kac-Moody algebra $\frak{g}$. In \cite{Kas94'},
Kashiwara introduced a class of integrable modules over
$U_q(\frak{g})$ called extremal weight modules, which is a natural
generalization of integrable highest weight or lowest weight modules. There exist
a global crystal base and a crystal base of an extremal weight
module, and the crystal base of an extremal weight crystal, simply
an extremal weight crystal, appears as a subcrystal of that of the
modified quantized enveloping algebra $\widetilde{U}_q(\frak{g})$
\cite{Kas94'}. Suppose that $\frak{g}$ is an affine Kac-Moody
algebra of finite rank. Then an extremal weight crystal of positive
(resp. negative) level is isomorphic to the crystal base of an
integrable highest (resp. lowest) weight module. In \cite{AK,
Kas02}, a level zero extremal weight module has been
studied in detail, and it was conjectured by Kashiwara \cite{Kas02}
that the structure of a level zero  extremal weight crystal can be described in terms of  a tensor product of the crystal bases of level zero fundamental
weight modules and Laurent Schur polynomials. In
\cite{BN}, Beck and Nakajima proved this conjecture (see also
\cite{B,Na} for the case when $\frak{g}$ is symmetric), and
furthermore based on the study of extremal weight modules they
proved the Kashiwara's conjecture on Peter-Weyl decomposition and
the Lusztig's conjecture on two-sided cell structures of
$\widetilde{U}_q(\frak{g})$ in a purely algebraic way, though there
is a geometric background related with quiver varieties.

A natural question arises whether there is also a nice description of
extremal weight crystals and their tensor products when $\frak{g}$
is an infinite rank affine Lie algebra. In \cite{K09} the author
studied extremal weight crystals of type $A_{+\infty}$, and he showed
that an extremal weight crystal is isomorphic to the tensor product
of a lowest weight crystal and a highest weight crystal, and the
Grothendieck ring generated by the isomorphism classes of extremal
weight crystals is isomorphic to the Weyl algebra of infinite rank.
A Littlewood-Richardson rule of extremal weight crystals is then
described explicitly using the operators induced from the
multiplication by Schur functions together with their adjoints.

The purpose of this paper is to study extremal weight crystals of
type $A_\infty$ (or extremal weight $\gl_\infty$-crystals), where
$\gl_\infty$ denotes the infinite rank affine Lie algebra of type
$A_\infty$.
In this case, it should be noted that (1) an extremal weight crystal
is always connected (Proposition \ref{connectedness}), (2) there are extremal weight crystals of
non-zero level, which are not isomorphic to a highest or lowest
weight crystal (Theorem \ref{realization thm}), and (3) the tensor product of two extremal weight crystals of non-negative level is a disjoint union of extremal weight crystals of non-negative level (Theorem \ref{category C}). These are
important features of extremal weight $\gl_\infty$-crystals, which
do not necessarily hold in the affine types of finite rank. Also,
as in the case of $A_{+\infty}$, we need certain non-commuting operators
to describe the tensor product of extremal $\gl_\infty$-crystals
because of the non-existence of characters and the non-commutativity
of tensor products.

Let us explain our results in detail. For an integral weight
$\Lambda$ of level $k\geq 0$, we denote by $\B(\Lambda)$ the crystal base of the
extremal weight module over $U_q(\gl_\infty)$ with extremal weight $\Lambda$. Let $\B$ be the crystal base of the natural representation of
$U_q(\frak{gl}_\infty)$. The connected components of $\B^{\otimes
n}$ ($n\geq 1$) are parameterized by partitions $\lambda$ of $n$,
say $\B_{\lambda}$. Note that the crystal $\B_\lambda$ is not
isomorphic to a highest weight or lowest weight
crystal. Then we show that $\B(\Lambda)$ is connected and there exist
unique partitions $\mu,\nu$ and the dominant integral weight
$\Lambda'$ of level $k$ such that
\begin{equation*}
\B(\Lambda)\simeq \B_\mu\otimes\B_\nu^\vee\otimes \B(\Lambda'),
\end{equation*}
where $\B_\nu^\vee$ is the dual crystal of $\B_\nu$ (Theorem
\ref{realization thm}). Note that $\B(\Lambda')$ is the crystal base of
the integrable highest weight module with highest weight $\Lambda'$. An
extremal weight crystal of non-positive level is also characterized
from the above result by taking its dual.

Next, we consider a category $\mathcal{C}$ of
$\frak{gl}_\infty$-crystals, whose objects are disjoint unions of
extremal weight crystals of non-negative level  with certain finite conditions on the multiplicity of connected components (see Section \ref{Cat C}). We show that
$\mathcal{C}$ is a monoidal category under tensor product of
crystals (Theorem \ref{category C}). We remark that
the tensor product in $\mathcal{C}$ is not necessarily commutative, that is, there is not always an isomorphism between $B\otimes B'$ and $B'\otimes B$ for any two objects $B$, $B'$ in $\mathcal{C}$ (see for example, Proposition \ref{Pieri}).

The Grothendieck group $\K$ of the category $\mathcal{C}$ admits a natural structure of a non-commutative associative $\Z$-algebra with $1$ induced from the monoidal structure of $\mathcal{C}$.
 Let $z=\{\,z_k\,|\,k\in\Z\,\}$ be a set of formal
commuting variables, and let $\mathscr{R}$ be the ring of formal
power series in $z$ with coefficients in $\mathbb{Z}$.   Now, let
$\cD$ be an Ore extension of $\R_\mathbb{Q}=\mathbb{Q}\otimes_{\mathbb{Z}} \mathscr{R}$ associated with a commuting family
of derivations $\gamma^\pm_n=(-1)^{n-1}\sum_{k\in\Z}z_{k\mp
n}\frac{\partial}{\partial z_k}$ ($n\in\Z$). We may view
$\mathscr{D}$ as a non-commutative polynomial ring over
$\mathscr{R}_\mathbb{Q}$ in $\gamma^{\pm}=\{\,\gamma^\pm_n\,|\,n\in\Z\,\}$.
Then we show that there exists a $\mathbb{Q}$-algebra isomorphism
$$ \K_\mathbb{Q} \stackrel{\sim}{\longrightarrow} \mathscr{D}^{\rm \, opp}$$ (Theorem
\ref{K-main}), where isomorphism classes of an integrable highest
weight crystal and a level zero extremal weight crystal are mapped
to polynomials in $z$ and $\gamma^\pm$, respectively. Here $\K_\mathbb{Q} =\mathbb{Q}\otimes_{\mathbb{Z}}\K$ and
$\mathscr{D}^{\rm \, opp}$ denotes the opposite algebra of
$\mathscr{D}$.

Based on the above results, we obtain a Littlewood-Richardson rule
for extremal weight crystals of non-negative level (Theorem
\ref{Extremal LR rule}), which is given explicitly in terms of
classical Littlewood-Richardson coefficients. In fact, the tensor
product of level zero extremal weight crystals corresponds to the
product of double symmetric functions, whose decomposition can be
given by the classical Littlewood-Richardson rule due to a crystal
$(\gl_\infty,\gl_n)$-duality on the $n$-th exterior algebra generated
by the natural representation of $\gl_\infty$ (Proposition
\ref{duality on En}), and the decomposition of the tensor product of
integrable highest weight crystals is explained by
using a crystal version of the $(\gl_\infty,\gl_n)$-duality  on the
level $n$ fermionic Fock space due to Frenkel \cite{Fr} (Proposition
\ref{duality-1}) following \cite{W99'}. Hence the only non-trivial
part is the decomposition of the tensor product
$\B(\Lambda)\otimes\B_{\mu}\otimes\B_{\nu}^\vee$, where $\Lambda$ is
a dominant integral weight and $\mu,\nu$ are partitions, and it is
obtained by analyzing commutation relations for
monomials in $z$ and $\gamma^\pm$, equivalently Pieri rules for
extremal weight crystals (Proposition \ref{Pieri}).

Finally, we discuss some applications. Let $\mathcal{C}^\vee$ be the category of
$\gl_\infty$-crystals consisting of dual crystals $B^\vee$ for
$B\in\mathcal{C}$ and let $\mathcal{C}^{\rm l.w.}$ be a full
subcategory of $\mathcal{C}^\vee$ whose objects are disjoint union of integrable  lowest weight crystal.
We denote by $\K^\vee$ and
$\K^{\rm l.w.}$ the corresponding Grothendieck groups. Then we consider
 a left $\K^\vee$-module structure on $\K^{\rm
l.w.}$, which is induced from the action of $\mathscr{D}$ on $\mathscr{R}_{\mathbb{Q}}$ as
differential operators, or from the composite of following two functors;
$$\mathcal{C}^\vee \times \mathcal{C}^{\rm l.w.} \stackrel{\otimes}{\longrightarrow}
\mathcal{C}^\vee \stackrel{{\rm pr}}{\longrightarrow}
\mathcal{C}^{\rm l.w.},$$ where ${\rm pr}$ is the natural
projection functor. Using the Littlewood-Richardson rule, we obtain an
explicit combinatorial description of the action of $\K^\vee$ on
$\K^{\rm l.w.}$. We observe that the action of level zero extremal
weight crystals are transitive on the set of integrable lowest
weight crystals with a fixed level. As an application, we obtain a
new interpretation of a one-to-one correspondence between level $n$
integrable highest (or lowest) weight $\gl_\infty$-modules and
finite dimensional $\gl_n$-modules which comes from the
$(\gl_\infty,\gl_n)$-duality on level $n$ fermionic Fock space \cite{Fr}
(Theorem \ref{K0 module structure}). As another application, we
construct an action of the Hall-Littlewood vertex operators \cite{J}
on $\Z[q]\otimes_{\mathbb{Z}} \K^{\rm l.w.}$ (Theorem \ref{Vertex operator}),
which naturally yields an $A_\infty$-analogue of Hall-Littlewood
function.


For a combinatorial realization, most crystals in this
paper are embedded in a set of binary matrices of various shapes,
equivalently an (infinite) abacus model. Also, two kinds of Kashiwara operators on binary matrices \cite{DK,K08-JACO,La}, which produces various dualities,
play a crucial role in proving our main results, while the rational
semistandard tableaux for $\gl_n$ \cite{St,Str} were used to
understand extremal weight crystals of type $A_{+\infty}$ \cite{K09}.

We also expect a similar result for the other infinite rank affine
Lie algebras, that is, roughly speaking, the Grothendieck ring
generated by extremal weight crystals can be realized as a ring of
differential operators acting on the character ring of integrable highest weight or lowest weight modules.  

The paper is organized as follows. In Section 2, we review briefly
the notion of crystals. In Section 3, we introduce a double crystal
(or bicrystal) structure on binary matrices, which is our main
method. In Section 4, we give a characterization of extremal weight
crystals. In Section 5, we introduce the monoidal category
$\mathcal{C}$, characterize its Grothendieck ring, and give a
Littlewood-Richardson rule for extremal weight crystals. In Section
6, we study the action of $\K^\vee$ on $\K^{\rm l.w.}$  and discuss its applications.\vskip 2mm

{\bf Acknowledgement} The author would like to thank the referee for careful reading of the first manuscript and many helpful comments on it.

\section{Crystals}
\subsection{Review on crystals}\label{crystal}
Let $I$ be an index set. Let $\mathfrak{g}$ be a symmetrizable Kac-Moody algebra associated with a generalized Cartan matrix $A=(a_{ij})_{i,j\in I}$. Denote the weight lattice of $\mathfrak{g}$ by $P$, the set of simple roots by $\Pi=\{\,\alpha_i\,|\,i\in I\,\}\subset P$, and the set of simple coroots by $\Pi^\vee=\{\,h_i\,|\,i\in I\,\}\subset P^\vee$ with $\langle\alpha_j,h_i\rangle=a_{ij}$.

Let $U_q(\frak{g})$ be the  quantized enveloping algebra associated with $\frak{g}$ introduced by Drinfeld and Jimbo. In \cite{Kas91}, Kashiwara introduced the notion of crystal base of a $U_q(\frak{g})$-module $V$, which can be can be viewed as a limit of $V$ at $q=0$. The crystal base is an $I$-colored oriented graph with important combinatorial information of $V$. The existence of the crystal bases of $U_q(\frak{g})$-modules which are related with the work in this paper has been proved in \cite{Kas91, Kas93, Kas94',KN}.

Based on properties of crystal bases, one can define the notion of $\frak{g}$-crystal (or crystal for short) as follows  (see \cite{Kas94} for a general review and references therein).

A {\it $\mathfrak{g}$-crystal} is a set
$B$ together with the maps ${\rm wt} : B \rightarrow P$,
$\varepsilon_i, \varphi_i: B \rightarrow \mathbb{Z}\cup\{-\infty\}$ and
$\te_i, \tf_i: B \rightarrow B\cup\{{\bf 0}\}$ ($i\in I$) such that
for $b\in B$ and $i\in I$
\begin{itemize}
\item[(1)] $\varphi_i(b) =\langle {\rm wt}(b),h_i \rangle +
\varepsilon_i(b),$

\item[(2)]  ${\rm wt}({\te_i}b)={\rm wt}(b)+\alpha_i$ if
$\te_i b \neq {\bf 0}$, and ${\rm wt}({\tf_i}b)={\rm
wt}(b)-\alpha_i$ if $\tf_i b \neq {\bf 0}$,

\item[(3)] $\tf_i b = b'$ if and only if $b = \te_i
b'$ for $b, b' \in B$,
\end{itemize}
where ${\bf 0}$ is a formal symbol. Here we assume that $-\infty+n=-\infty$ for all $n\in\Z$. Note that $B$
is equipped with an $I$-colored oriented graph structure, where
$b\stackrel{i}{\rightarrow}b'$ if and only if $b'=\tf_{i}b$ for
$b,b'\in B$ and $i\in I$. We call $B$ {\it connected} if it is connected as a graph.
\vskip 2mm

Let us review some terminologies on crystals.
Let $B$ be a crystal.  We say that $B$ is {\it normal} if $\varepsilon_i(b)={\rm
max}\{\,k\,|\,\te_i^kb\neq {\bf 0}\,\}$ and $\varphi_i(b)={\rm
max}\{\,k\,|\,\tf_i^kb\neq {\bf 0}\,\}$ for $b\in B$ and $i\in I$, and put $\te_i^{\rm max}b =
\te_i^{\varepsilon_i(b)}b$ and $\tf_i^{\rm max}b =\tf_i^{\varphi_i(b)}b$.  We say that $B$ is {\it regular} if  $B$ is as a $\g_J$-crystal, isomorphic to the crystal base of an integrable $U_q(\g_J)$-module for $J\subset I$ such that $\{\,\alpha_i\,|\,i\in J\,\}$ is of finite type, where $\g_J$ is the Kac-Moody algebra associated with $A_J=(a_{ij})_{i,j\in J}$ \footnote{Erratum : the definition of regular crystal in \cite[p1343]{K09}, which was given as that of normal crystal,  should be replaced with the one given here.}.

Let $W$ be the Weyl group of $\g$, that is, the subgroup of $GL(P)$ generated
by $r_i$ ($i\in I$), where $r_i$ is the simple reflection given by
$r_i(\lambda)=\lambda-\langle \lambda,h_i \rangle\alpha_i$ for
$\lambda\in P$.
 A regular crystal $B$ admits an action of the Weyl group on $B$ as follows; for
$i\in I$ and $b\in B$
\begin{equation}\label{Weyl}
S_{r_i} b=
\begin{cases}
\tf_i^{\,\langle {\rm wt} (b),h_i \rangle}b, & \text{if $\langle {\rm
wt} (b),h_i \rangle\geq 0$}, \\
\te_i^{\,-\langle {\rm wt} (b),h_i \rangle}b, & \text{if $\langle {\rm
wt} (b),h_i \rangle\leq 0$}, \\
\end{cases}
\end{equation}
and $S_w=S_{r_{i_1}}\cdots S_{r_{i_t}}$ for  $w\in W$ with a reduced expression $w=r_{i_1}\cdots r_{i_t}$.
Note that if $B$ is regular, then $B$ is normal.

The {\it dual crystal $B^\vee$ of $B$} is defined to be the set
$\{\,b^\vee\,|\,b\in B\,\}$ with{\allowdisplaybreaks
\begin{equation}
\begin{aligned}
&{\rm wt}(b^\vee)=-{\rm wt}(b), \\
&\varepsilon_i(b^\vee)=\varphi_i(b), \ \
\varphi_i(b^\vee)=\varepsilon_i(b), \\
&\te_i(b^\vee)=\left(\tf_i b \right)^\vee, \ \
\tf_i(b^\vee)=\left(\te_i b \right)^\vee,
\end{aligned}
\end{equation}}
for $b\in B$ and $i\in I$. We assume that ${\bf 0}^\vee={\bf
0}$.

Let $B_1$ and $B_2$ be crystals. A {\it morphism}
$\psi : B_1 \rightarrow B_2$ is a map from $B_1\cup\{{\bf 0}\}$ to
$B_2\cup\{{\bf 0}\}$ such that
\begin{itemize}
\item[(1)] $\psi(\bf{0})=\bf{0}$,

\item[(2)] ${\rm wt}(\psi(b))={\rm wt}(b)$,
$\varepsilon_i(\psi(b))=\varepsilon_i(b)$, and
$\varphi_i(\psi(b))=\varphi_i(b)$ if $\psi(b)\neq \bf{0}$,

\item[(3)] $\psi(\te_i b)=\te_i\psi(b)$ if $\psi(b)\neq \bf{0}$ and
$\psi(\te_i b)\neq \bf{0}$,

\item[(4)] $\psi(\tf_i
b)=\tf_i\psi(b)$ if $\psi(b)\neq
\bf{0}$ and $\psi(\tf_i b)\neq \bf{0}$,
\end{itemize}
for $b\in B_1$ and $i\in I$.
We call $\psi$ an {\it embedding} and $B_1$ a {\it subcrystal of}
$B_2$ when $\psi$ is injective, and call $\psi$ {\it strict} if
$\psi : B_1\cup\{{\bf 0}\} \rightarrow B_2\cup\{{\bf 0}\}$ commutes
with $\te_i$ and $\tf_i$ for all $i\in I$, where we assume that $\te_i{\bf
0}=\tf_i{\bf 0}={\bf 0}$.

A {\it tensor product $B_1\otimes B_2$}  is defined to be $B_1\times B_2$ as a set with elements  denoted by $b_1\otimes b_2$, where  {\allowdisplaybreaks
\begin{equation}
\begin{split}
{\rm wt}(b_1\otimes b_2)&={\rm wt}(b_1)+{\rm wt}(b_2), \\
\varepsilon_i(b_1\otimes b_2)&= {\rm
max}\{\varepsilon_i(b_1),\varepsilon_i(b_2)-\langle {\rm
wt}(b_1),h_i\rangle\}, \\
\varphi_i(b_1\otimes b_2)&= {\rm max}\{\varphi_i(b_1)+\langle {\rm
wt}(b_2),h_i\rangle,\varphi_i(b_2)\},\\
{\te}_i(b_1\otimes b_2)&=
\begin{cases}
{\te}_i b_1 \otimes b_2, & \text{if $\varphi_i(b_1)\geq \varepsilon_i(b_2)$}, \\
b_1\otimes {\te}_i b_2, & \text{if
$\varphi_i(b_1)<\varepsilon_i(b_2)$},
\end{cases}\\
{\tf}_i(b_1\otimes b_2)&=
\begin{cases}
{\tf}_i b_1 \otimes b_2, & \text{if  $\varphi_i(b_1)>\varepsilon_i(b_2)$}, \\
b_1\otimes {\tf}_i b_2, & \text{if $\varphi_i(b_1)\leq
\varepsilon_i(b_2)$},
\end{cases}
\end{split}
\end{equation}
\noindent for $i\in I$. Here we assume that ${\bf 0}\otimes
b_2=b_1\otimes {\bf 0}={\bf 0}$.} Then $B_1\otimes B_2$ is a
crystal. Note that $(B_1\otimes B_2)^\vee\simeq
B_2^\vee\otimes B_1^\vee$.

For $b_i\in B_i$ ($i=1,2$), we say that {\it $b_1$ is
{\rm ($\g$-)}equivalent to $b_2$}, and write $b_1 \equiv b_2$
if there exists a crystal isomorphism $C(b_1)\rightarrow
C(b_2)$ sending $b_1$ to $b_2$, where $C(b_i)$ denotes the connected
component of $B_i$ including $b_i$ ($i=1,2$).

\subsection{The Lie algebra $\gl_\infty$}
Let $\gl_\infty$ denote the Lie algebra of complex matrices
$(a_{ij})_{i,j\in \Z}$ with finitely many non-zero entries.

Let $E_{ij}$ be the elementary matrix with $1$ at the $i$-th row and
the $j$-th column and zero elsewhere.
Let $\h=\bigoplus_{i\in \Z}\C E_{ii}$ be the Cartan subalgebra of
$\gl_\infty$ and $\langle\cdot,\cdot\rangle$ denote the natural
pairing on $\frak{h}^*\times\frak{h}$. We denote by $\Pi^{\vee}=\{\,
h_i=E_{ii}-E_{i+1,i+1}\,|\, i\in\Z\, \}$ the set of simple coroots,
and $\Pi=\{\, \alpha_i=\epsilon_i-\epsilon_{i+1} \,|\, i\in\Z\, \}$
the set of simple roots, where $\epsilon_i\in\h^*$ is determined by
$\langle \epsilon_i,E_{jj}\rangle=\delta_{ij}$.

Let
$P=\bigoplus_{i\in\Z}\Z\epsilon_i\oplus\Z\Lambda_0=\bigoplus_{i\in\Z}\Z\Lambda_i$ be
the weight lattice of $\gl_{\infty}$, where $\Lambda_0$ is defined
by
$\langle\Lambda_0,E_{-j+1,-j+1}\rangle=-\langle\Lambda_0,E_{jj}\rangle=\frac{1}{2}$
for $j\geq 1$, and $\Lambda_i$ is given by
$\Lambda_0-\sum_{k=i}^{0}\epsilon_k$ (resp.
$\Lambda_0+\sum_{k=1}^i\epsilon_k$) for $i<0$ (resp. $i>0$). We call
$\Lambda_i$ the $i$-th fundamental weight. A partial ordering on $P$
is defined as usual.

For $k\in\Z$, let $P_k=k\Lambda_0+\bigoplus_{i\in\Z}\Z\epsilon_i$ be
the set of integral weights of level $k$. Let $P^+=\{\,\Lambda\in
P\,|\,\langle\Lambda,h_i \rangle\geq 0,\
i\in\Z\,\}=\sum_{i\in\Z}\Z_{\geq 0}\Lambda_i$ be the set of dominant
integral weights. We also put $P^+_k=P^+\cap P_k$ for $k\geq 0$ (note that $P^+_0=\{0\}$).
Note that for $\Lambda=\sum_{i\in\Z}c_i\Lambda_i\in P$, the level of
$\Lambda$ is $\sum_{i\in\Z}c_i$ since
$\epsilon_i=\Lambda_i-\Lambda_{i-1}$ for $i\in\Z$. If we put
$\Lambda_\pm=\sum_{i; c_i\gtrless 0}c_i\Lambda_i$, then
$\Lambda=\Lambda_+-\Lambda_-$ with $\Lambda_\pm\in P^+$.

For $n\geq 1$, let
$\Z_+^n=\{\,(\lambda_1,\cdots,\lambda_n)\,|\,\lambda_i\in\Z,\
\lambda_1\geq \cdots\geq \lambda_n\,\}$ be the set of {generalized
partitions of length $n$}. For $\lambda\in\Z_+^n$, we put
\begin{equation}
\begin{split}
\Lambda_{\lambda}& =\Lambda_{\lambda_1}+\cdots+\Lambda_{\lambda_n}.
\end{split}
\end{equation}
Note that $\Z_+^n$ parameterizes $P^+_n$.

Let $I$ be an interval in $\mathbb{Z}$. We denote by $\gl_I$ the
subalgebra of $\gl_{\infty}$ spanned by $E_{ij}$ for $i,j\in I$.  For
$p,q\in\Z$, we put $[p,q]=\{\,p,p+1,\ldots,q\,\}$ ($p<q$),
$[p,\infty)=\{\,p,p+1,\ldots\}$ and
$(-\infty,q\,]=\{\,\ldots,q-1,q\,\}$. For $n\geq 1$, we denote
$[1,n]$ by $[n]$ for simplicity.

For $\Lambda\in P^+$, we denote by $\B(\pm\Lambda)$ the crystal base
of the irreducible $U_q(\gl_{\infty})$-module with highest (resp.
lowest) weight vector $u_{\pm\Lambda}$ of weight $\pm\Lambda$, which
are connected regular $\gl_{\infty}$-crystals. We denote by $\B$ and $\B^\vee$ the crystal base of the
natural representation of $U_q(\gl_{\infty})$ and its dual
respectively, which are also connected regular crystals. The associated colored oriented graphs are
{\allowdisplaybreaks
\begin{align*}\label{natural representations}
\B: \ \ &\ \ \ \, \cdots \stackrel{-3}{\longrightarrow}
-2\stackrel{-2}{\longrightarrow}-1\stackrel{-1}{\longrightarrow} 0
\stackrel{0}{\longrightarrow}1 \stackrel{1}{\longrightarrow} 2
\stackrel{2}{\longrightarrow}3\stackrel{3}{\longrightarrow}
\cdots,\\
\B^\vee:\ \ & \cdots \stackrel{3}{\longrightarrow}
3^\vee\stackrel{2}{\longrightarrow}2^\vee\stackrel{1}{\longrightarrow}
1^\vee \stackrel{0}{\longrightarrow}
  0^\vee \stackrel{-1}{\longrightarrow}
-1^\vee\stackrel{-2}{\longrightarrow}-2^\vee\stackrel{-3}{\longrightarrow}
\cdots,
\end{align*}}
where ${\rm wt}(i)=\epsilon_i$ for $i\in\Z$.

\section{Double crystal structure on binary matrices}

\subsection{Crystal operators on binary matrices}\label{matrices}

For intervals $I,J$ in $\Z$, let $\bM_{I,J}$ be the set of $I\times
J$ matrices $A=(a_{ij})$ with $a_{ij}\in\{\,0,1\,\}$. We denote by
$A_i$ the $i$-th row of $A$ for
$i\in I$.

Suppose that $|J|\geq 2$ and let $J^\circ=\{\,k\,|\,k,k+1\in J\,\}$.
For $i\in I$ and $k\in J^\circ$, let
\begin{equation*}
\sigma_k(A_i)=
\begin{cases}
 \ + &\  \text{if $(a_{i\, k},a_{i\, k+1})=(1,0)$}, \\
 \ - &\ \text{if $(a_{i\, k},a_{i\, k+1})=(0,1)$}, \\
 \ \ \cdot &\  \text{otherwise}.
\end{cases}
\end{equation*}
We say that $A$ is {\it  row $k$-admissible} if there exist $M, N\in I$ ($M\leq N$) such that $\sigma_k(A_i)\neq +$ for all $i<M$ and  $\sigma_k(A_i)\neq -$ for all $i>N$. Note that if $I$ is finite, then $A$ is row $k$-admissible.

Let us define the operators $\te_k, \tf_k$ for $k\in J^\circ$ on the set of row $k$-admissible matrices in $\bM_{I,J}$.
Suppose that $|I|=1$. For $A\in \bM_{I,J}$, we define
\begin{equation}
\begin{split}
\te_kA&=
\begin{cases}
A+E_k-E_{k+1}, & \text{if $\sigma_k(A)=-$}, \\
{\bf 0},  & \text{otherwise},
\end{cases} \\
\tf_kA&=
\begin{cases}
A-E_k+E_{k+1}, & \text{if $\sigma_k(A)=+$}, \\
{\bf 0},  & \text{otherwise},
\end{cases}
\end{split}
\end{equation}
where $E_k$ is the $I\times J$ matrix with $1$ in the $k$-th column and $0$ elsewhere.

Suppose that $|I|\geq 2$. For a row $k$-admissible $A\in \bM_{I,J}$, consider the sequence
$$\sigma_k(A)=(\sigma_k(A_i))_{i\in I}=(\ldots,\sigma_k(A_{i-1}),\sigma_k(A_i),\sigma_k(A_{i+1}),\ldots).$$ We replace a pair
$(\sigma_k(A_s),\sigma_k(A_{s'}))=(+,-)$ such that $s<s'$ and $\sigma_k(A_{t})=\cdot$ for $s<t<s'$ by
$(\,\cdot\,,\,\cdot\,)$ in $\sigma_k(A)$, and repeat this process
as far as possible until we get a sequence $\td{\sigma}_k(A)=(\td{\sigma}_k(A_i))_{i\in I}$ with no $+$ placed to the
left of $-$. It is not difficult to see that this procedure must end after a finite number of steps since $A$ is row $k$-admissible.

Then we define $\te_k A$ (resp. $\tf_k
A$) to be the matrix in $\bM_{I,J}$ given by applying $\te_k$ (resp.
$\tf_k$) to the row of $A$ corresponding to the right-most $-$  (resp. the left-most $+$) in $\td{\sigma}_k(A)$, and
$\te_kA={\bf 0}$ (resp. $\tf_kA={\bf 0}$) if there is no such row. Note that $\te_kA$ (resp. $\tf_kA$) is also row $k$-admissible if $\te_kA\neq {\bf 0}$ (resp. $\tf_kA\neq {\bf 0}$).

\begin{rem}{\rm
We may regard each $A_i$ as an element of the crystal base of the fundamental weight representation or the trivial representation of $U_q(\gl_{\{k,k+1\}})$. Hence if $I$ is a finite set, then the actions of $\te_k$ and $\tf_k$ on $A$ defined here coincides with those on the tensor product of $A_i$'s, and our definition is equal to the algorithm usually known as {\it signature rule} (cf. \cite{KN}). But, we do not always have a crystal structure on the set of row $k$-admissible matrices, since the weight of $A$ with respect to $\gl_{\{k,k+1\}}$ is not always defined naturally when $I$ is an infinite set.
}
\end{rem}

Let
$$\rho : \bM_{I,J} \longrightarrow \bM_{-J,I}$$ be a bijection given by
$\rho((a_{ij}))=(a'_{-j\,i})\in \bM_{-J,I}$ with $a'_{-j\,i}=a_{ij}$, where $-J=\{\,-j\,|\,j\in
J\,\}$.

Suppose that $|I|\geq 2$.
For $A\in \bM_{I,J}$ and $l\in I^\circ$, we say that $A$ is {\it column $l$-admissible} if $\rho(A)$ is row
$l$-admissible. Then for a column $l$-admissible $A\in\bM_{I,J}$,
we define
\begin{equation}
\begin{split}
\tE_l(A)=\rho^{-1}\left(\te_{l}\,\rho(A) \right), \ \ \tF_l(A)=\rho^{-1}\left(\tf_{l}\, \rho(A) \right)
\end{split}.
\end{equation}

The following lemma follows from \cite[Proposition 4.2]{K07} (see
also \cite{DK,La}, where essentially the same facts are stated in a
slightly different way).
\begin{lem}\label{commutativity}
Let $A\in \bM_{I,J}$ be given. If $A$ is both row $k$-admissible and column $l$-admissible for $l\in I^\circ$ and $k\in J^\circ$, then
$$\widetilde{x}_k\widetilde{X}_lA=\widetilde{X}_l\widetilde{x}_kA,$$
where $x=e,f$ and $X=E,F$.
\end{lem}

For $A\in\bM_{I,J}$, we say that {\it $A$ is row admissible} (resp.
{\it column admissible}) if $A$ is row $k$-admissible (resp. column
$l$-admissible) for all $k\in J^\circ$ (resp. $l\in I^\circ$). Note that if $I$ (resp. $J$) is a finite set, then
$A$ is always row (resp. column) admissible. Lemma \ref{commutativity} implies the following
immediately.

\begin{lem}\label{admissibility}
Let $A\in \bM_{I,J}$ be given. Suppose that $A$ is row admissible
and column $l$-admissible for $l\in I^\circ$. If $\widetilde{X}_lA\neq {\bf
0}$ $(X=E,F)$, then
\begin{equation*}
\td{x}_{k_1}\ldots\td{x}_{k_r}A \neq {\bf 0} \Longleftrightarrow \td{x}_{k_1}\ldots\td{x}_{k_r}(\td{X}_lA) \neq {\bf 0}
\end{equation*}
for $r\geq 1$ and $k_1,\ldots, k_r\in J^\circ$, where $x=e$ or $f$ for each $k_i$.
\end{lem}
We have a similar statement when $A$ is column admissible and
row $k$-admissible for $k\in J^\circ$.

For $A=(a_{ij})\in\bM_{I,J}$, put
\begin{equation}\label{Avee}
A^\vee=(1-a_{-i\,j})_{-i\in I, j\in J}\in \bM_{-I,J}.
\end{equation}
Suppose that $A$ and $A^\vee$ are row $k$-admissible for some $k\in J^\circ$.
Then we have
\begin{equation}\label{dual matrix}
\tf_k A = \left(\te_k A^\vee \right)^\vee, \ \ \ \ \te_k A = \left(\tf_k A^\vee \right)^\vee.
\end{equation}
The same statement holds for $\tE_l$ and $\tF_l$, when $A$ and
$A^\vee$ are column $l$-admissible for some $l\in I^\circ$. We call $A^\vee$ the
{\it dual of $A$}.

\subsection{Crystals of semistandard tableaux}
Let $\cP$ denote the set of partitions. We identify a partition
$\lambda=(\lambda_i)_{i\geq 1}$ with a Young diagram  as usual (see
\cite{Mac95}). The number of non-zero parts in $\lambda$ is denoted
by $\ell(\lambda)$. We denote by
$\lambda'$ the conjugate partition of $\lambda$. For a skew Young
diagram $\lambda/\mu$, $|\lambda/\mu|$ denotes the number of dots or
boxes in the diagram. Let $\A$ be a linearly ordered set. A tableau
$T$ obtained by filling $\lambda/\mu$ with entries in $\A$ is called
a {\it semistandard tableau of shape $\lambda/\mu$} if the entries
in each row are weakly increasing from left to right, and the
entries in each column are strictly increasing from top to bottom.
We denote by $SST_\A(\lambda/\mu)$ the set of all semistandard
tableaux of shape $\lambda/\mu$ with entries in $\A$. For $T\in
SST_{\A}(\lambda/\mu)$, let $w(T)_{\rm col}$ (resp. $w(T)_{\rm
row}$) denote the word of $T$ with respect to column (resp. row)
reading, where we read the entries column by column (resp. row by
row) from right to left (resp. top to bottom), and in each column
(resp. row)
 from top to bottom (resp. right to left).

Let $\A$ denote one of the crystals $\B$ and $\B^\vee$, which are
linearly ordered with respect to the partial ordering on $P$. Note
that the set of all finite words with letters in $\A$ is a
$\frak{gl}_{\infty}$-crystal, where each word of length $r\geq 1$ is
identified with an element in $\A^{\otimes
r}=\A\otimes\cdots\otimes\A$ ($r$ times). Given a skew Young diagram
$\lambda/\mu$, the injective image of $SST_\A(\lambda/\mu)$ in the
set of all finite words under the map $T \mapsto w(T)_{\rm col}$ (or
$w(T)_{\rm row}$) together with $\{{\bf 0}\}$ is invariant under
$\te_i,\tf_i$ ($i\in\Z$). Hence it is a regular $\frak{gl}_{\infty}$-crystal \cite{KN}.
In particular, for $\lambda\in\cP$, we have
$$SST_{\B}(\lambda)^\vee\simeq SST_{\B^\vee}(\lambda^\vee),$$
where $\lambda^\vee$ is the skew Young diagram obtained from
$\lambda$ by $180^\circ$-rotation.

For $\mu\in\cP$, we put
\begin{equation}
\cB_{\mu}=SST_{\B}(\mu),
\end{equation}
and we identify $\cB_\mu^\vee$ with $SST_{\B^\vee}(\mu^\vee)$. Note
that $\cB_{\mu}$ does not have  a highest weight or lowest weight element.

\begin{prop}\label{connectedness of bmn}  For $\mu,\nu\in\cP$, $\cB_{\mu}\otimes\cB_{\nu}^\vee$  and $\cB_{\nu}^\vee\otimes\cB_{\mu}$ are connected.
\end{prop}
\pf First, we claim that $\cB_\mu$ is connected. Suppose that $S,T\in\cB_\mu$ are given.
Choose $p\in\Z$ such that all entries in $S$ and $T$ are greater than $p$. Then $S$ is an element in $SST_{[p,\infty)}(\mu)$, which is a connected $\gl_{[p,\infty)}$-crystal with a unique highest weight element, say $u^{[p,\infty)}_{\mu}$ (see \cite{KN}). This implies that $S$ and $T$ are contained in the same connected component, and hence $\cB_\mu$ is connected.

Let $S\otimes T\in\cB_{\mu}\otimes\cB_{\nu}^\vee$ be given. Choose $p\in\Z$ such that $S\in SST_{[p,\infty)}(\mu)$. Then we have $\te_{i_1}\cdots\te_{i_r}S=u^{[p,\infty)}_{\mu}$  for some $i_1,\ldots,i_r\in [p,\infty)$.
By tensor product rule of crystals, we also have $$\te_{i_1}^{m_1}\cdots\te_{i_r}^{m_r}(S\otimes T)=u^{[p,\infty)}_{\mu}\otimes T'$$
for some $m_1,\ldots,m_r\geq 1$ and $T'\in \cB_{\nu}^\vee$.

If $p$ is sufficiently small, then we may assume that all the entries in $T$ (and hence in $T'$) are smaller than $(p+\ell(\mu))^\vee$. Choose $q$ such that $T'\in  SST_{(-\infty,q]^\vee}(\nu^\vee)$.
Note that $SST_{(-\infty,q]^\vee}(\nu^\vee)$ is a $\gl_{(-\infty,q]}$-crystal with a unique highest weight element $v^{(-\infty,q]}_{\nu}$.
Hence $\te_{j_1}\cdots\te_{j_s}T'=v^{(-\infty,q]}_{\nu}$ for some $j_1,\ldots,j_s\in (-\infty,q]$. Since $\{\,j_1,\ldots,j_s\,\}$ does not intersect with the entries in $u^{[p,\infty)}_{\mu}$, we have
$$\te_{j_1}\cdots\te_{j_s}\left(u^{[p,\infty)}_{\mu}\otimes T'\right)=u^{[p,\infty)}_{\mu}\otimes v^{(-\infty,q]}_{\nu}.$$

Now, let $U\otimes V\in \cB_{\mu}\otimes\cB_{\nu}^\vee$ be given. Then if $p$ is sufficiently small and $q$ is sufficiently large, then it follows from the same argument that $U\otimes V$ is also connected to $u^{[p,\infty)}_{\mu}\otimes v^{(-\infty,q]}_{\nu}$. This implies that $\cB_{\mu}\otimes\cB_{\nu}^\vee$ is connected.

The proof of the connectedness of $\cB_{\nu}^\vee\otimes\cB_{\mu}$ is almost identical and is omitted.
\qed\vskip 3mm

\begin{rem}\label{remark on bmn}{\rm
In the proof of Proposition \ref{connectedness of bmn}, we showed
that for $S\otimes T\in\B_{\mu}\otimes\B_{\nu}^\vee$ there exist
$k_1,\ldots,k_t\in\Z$ such that
$\te_{k_1}\cdots\te_{k_t}\left(S\otimes
T\right)=u^{[p,\infty)}_{\mu}\otimes v^{(-\infty,q]}_{\nu}$ for some
$p<q$. By applying suitable $\te_k$'s, we may assume that $p\ll 0\ll
q$.
 }
\end{rem}

Let $\E$ be the subset of $\bM_{\{1\},\Z}$ consisting of
$(a_{1j})_{j\in\mathbb{Z}}$ such that $\sum_{j}a_{1j}<\infty$. If we
define ${\rm wt}(A)=\sum_{j}a_{1j}\epsilon_j$ for
$A=(a_{1j})_{j\in\mathbb{Z}}\in\E$, then $\E$ is a regular
$\gl_{\infty}$-crystal with respect to $\te_k, \tf_k$ ($k\in \Z$) and ${\rm wt}$  (see Section \ref{matrices} for the
definitions of $\te_k$ and $\tf_k$ on $\bM_{\{1\},\Z}$). For $a\geq 1$, define
\begin{equation}\label{sigma}
\sigma_a : \cB_{(1^a)} \longrightarrow \E
\end{equation}
by $\sigma_a(S)=(a_{1j})_{j\in\Z}$ with ${\rm
wt}(S)=\sum_{j}a_{1j}\epsilon_j$.  It is easy to check  that $\sigma_a$ is a strict embedding and
\begin{equation}\label{E1 decomposition}
\E\simeq \bigsqcup_{a\geq 0}\cB_{(1^a)}.
\end{equation}
Indeed, $\E$ is the crystal base of the $q$-deformed exterior algebra of the natural representation of $U_q(\gl_\infty)$.

Let $\E^\vee=\{\,A^\vee\,|\,A\in \E\,\}$, where $A^\vee$ denotes the dual matrix of $A$  (\ref{Avee}).  If we
define ${\rm wt}(A^\vee)=\sum_{j}(a_{1j}-1)\epsilon_j$ for
$A=(a_{1j})_{j\in\mathbb{Z}}\in\E$, then $\E^\vee$ is a regular
$\gl_{\infty}$-crystal with respect to $\te_k, \tf_k$ ($k\in \Z$) and ${\rm wt}$, which is isomorphic to the dual crystal of $\E$.
Similarly,
for $b\geq 1$  we define
\begin{equation}
\tau_b : \cB_{(1^b)}^\vee \longrightarrow \E^\vee
\end{equation}
by $\tau_b(T)=(a_{1j})_{j\in\Z}$ with ${\rm
wt}(T)=\sum_{j}(a_{1j}-1)\epsilon_j$. Then $\tau_b$ is a strict embedding. For convenience, we assume that
$\sigma_0$ is a map sending trivial crystal to zero matrix in
$\bM_{\{1\},\Z}$, and $\tau_0$ is a map sending trivial crystal to
the matrix $(a_{1j})_{j\in\Z}$ with $a_{1j}=1$ for all $j\in\Z$.

\begin{lem}\label{babb commutes} For $a,b\geq 0$, we have
$$\cB_{(1^a)}\otimes \cB_{(1^b)}^\vee \simeq \cB_{(1^b)}^\vee\otimes \cB_{(1^a)}.$$
\end{lem}
\pf Consider $u^{[p,\infty)}_{(1^a)}\otimes v^{(-\infty,q]}_{(1^b)}$ (for
simplicity write $u^p_a\otimes v^q_b$) for some $p\ll 0 \ll q$. We
claim that $$u^p_a\otimes v^q_b\equiv v^q_b\otimes u^p_a,$$ which
implies that $\cB_{(1^a)}\otimes \cB_{(1^b)}^\vee \simeq
\cB_{(1^b)}^\vee\otimes \cB_{(1^a)}$ by Proposition
\ref{connectedness of bmn}.

Define
\begin{equation}\label{sigma times tau}
\begin{split}
\sigma_a\times\tau_b &: \cB_{(1^a)}\otimes \cB_{(1^b)}^\vee \longrightarrow \E\otimes\E^\vee, \\
\tau_b\times \sigma_a &:  \cB_{(1^b)}^\vee \otimes\cB_{(1^a)}
\longrightarrow \E^\vee\otimes\E
\end{split}
\end{equation}
by $\sigma_a\times\tau_b(S\otimes T)=\sigma_a(S)\otimes \tau_b(T)$ and $\tau_b\times \sigma_a(T\otimes S)=\tau_b(T)\otimes \sigma_a(S)$ for $S\in \cB_{(1^a)}$ and $T\in \cB_{(1^b)}^\vee$. Then $\sigma_a\times\tau_b$ and $\tau_b\times \sigma_a$ are strict embeddings.  Here we assume $\E\otimes\E^\vee$ and $\E^\vee\otimes\E$ as subsets of $ \bM_{[2],\Z}$, where $A_1\otimes A_2$ is identified with a matrix $A \in \bM_{[2],\Z}$ whose $i$-th row  is $A_i$ ($i=1,2$).
Note that $(\sigma_a\times\tau_b)(u^p_a\otimes v^q_b)=A=(a_{ij})$,
where $a_{1j}=1$ if and only if $p\leq j\leq p+a-1$, and $a_{2j}=0$
if and only if $q-b+1\leq j\leq q$, while
$(\tau_b\times\sigma_a)(v^q_b\otimes u^p_a)=B=(b_{ij})$ with
$b_{ij}=a_{3-i\, j}$ for all $i,j$.

Choose $r\ll p$ and $s\gg q$. Let $\pi_{[r,s]}
: \bM_{[2],\Z} \rightarrow \bM_{[2],[r,s]}$ be the restriction map
sending a matrix to its $[2]\times [r,s]$ submatrix. Then $\pi_{[r,s]}(A)$ is column
$1$-admissible and $\tE_1^{\rm max} \pi_{[r,s]}(A)=\pi_{[r,s]}(B)$. By Lemma
\ref{admissibility}, we have
\begin{equation*}
\widetilde{x}_{i_1}\cdots\widetilde{x}_{i_t} \pi_{[r,s]}(A) \neq {\bf 0} \Longleftrightarrow \widetilde{x}_{i_1}\cdots\widetilde{x}_{i_t} \pi_{[r,s]}(B)\neq {\bf 0},
\end{equation*}
and hence
\begin{equation*}
\widetilde{x}_{i_1}\cdots\widetilde{x}_{i_t}A \neq {\bf 0} \Longleftrightarrow \widetilde{x}_{i_1}\cdots\widetilde{x}_{i_t} B\neq {\bf 0}
\end{equation*}
for $t\geq 1$ and $r\leq i_1,\ldots,i_t\leq s-1$, where $x=e$ or $f$ for each $i_k$.
Since $r$ and $s$ are arbitrary and ${\rm wt}(A)={\rm wt}(B)$, we have $A\equiv B$, which implies that $u^p_a\otimes v^q_b\equiv v^q_b\otimes u^p_a$. \qed\vskip 3mm

\begin{rem}\label{remark on babb commutes}{\rm
More generally, we can check by the argument in Lemma \ref{babb
commutes} that for $S\otimes T\in \cB_{(1^a)}\otimes
\cB_{(1^b)}^\vee$, if ${\rm
wt}(S)=\epsilon_{j_1}+\cdots+\epsilon_{j_a}$ and ${\rm
wt}(T)=-\epsilon_{j'_1}-\cdots-\epsilon_{j'_b}$ with
$j_1<\cdots<j_a<j'_1<\cdots<j'_b$, then $S\otimes T\equiv T\otimes
S$. }
\end{rem}

For $n\geq 1$, let $\E^n$ be the subset of $\bM_{[n],\Z}$ consisting
of matrices $A=(a_{ij})$ such that $A_i\in\E\subset \bM_{\{i\},\Z}$
for all $i\in [n]$. Then $\E^n$ is row admissible and can be identified with
$\E^{\otimes n}$ as a $\gl_\infty$-crystal, where $A$ is identified with $A_1\otimes\cdots\otimes A_n$. Also we may identify the dual crystal $\left(\E^n\right)^\vee$ with the set $\{\,A^\vee\,|\,A\in\E^n\,\}$, where $A^\vee$ denotes the dual matrix of $A$.

Let $\mu,\nu\in\cP$ be given with $\ell(\mu')=m$ and $\ell(\nu')=n$.
We may regard $\cB_{\mu}\subset
\cB_{(1^{\mu'_m})}\otimes\cdots\otimes \cB_{(1^{\mu'_1})}$, where
the $k$-th column of $S\in \cB_\mu$ (from the right) is an element
in $\cB_{(1^{\mu'_{m-k+1}})}$. Composing with (\ref{sigma}), we
have a strict embedding
\begin{equation}\label{sigma'}
\sigma_{\mu}=\sigma_{\mu'_m}\times\cdots\times\sigma_{\mu'_1} :
\cB_\mu \longrightarrow  \E^m
\end{equation}
(cf.(\ref{sigma times tau})). Similarly, we may regard $\cB_{\nu}^\vee\subset
\cB_{(1^{\nu_1'})}^\vee\otimes\cdots\otimes
\cB_{(1^{\nu_n'})}^\vee$, and have a strict embedding
\begin{equation}
\tau_{\nu}=\tau_{\nu'_1}\times\cdots\times\tau_{\nu'_n} : \cB_\nu^\vee
\longrightarrow \left(\E^n\right)^\vee.
\end{equation}

\begin{prop}\label{commmutativity of BmBnvee} For $\mu,\nu\in\cP$, we have
$$\cB_{\mu}\otimes\cB_{\nu}^\vee\simeq \cB_{\nu}^\vee\otimes \cB_{\mu}.$$
\end{prop}
\pf Consider  $u^{[p,\infty)}_{\mu}\otimes v^{(-\infty,q]}_{\nu}$
(see the proof of Proposition \ref{connectedness of bmn} for its
definition). We assume that $p\ll 0\ll q$. Let $A$ be the image of
$u^{[p,\infty)}_{\mu}\otimes v^{(-\infty,q]}_{\nu}$ under the strict embedding
$$\sigma_{\mu}\times \tau_{\nu} : \cB_{\mu}\otimes\cB_{\nu}^\vee \longrightarrow \E^m\otimes \left(\E^n\right)^\vee,$$
where $\ell(\mu')=m$ and $\ell(\nu')=n$. Let us write
$A=(A_1,\ldots,A_m,A_{m+1},\ldots,A_{m+n})\in \bM_{[m+n],\Z}$. Note that $A_i$ ($1\leq
i\leq m$) corresponds to the $i$-th column of $u^{[p,\infty)}_{\mu}$
and $A_{m+j}$ ($1\leq j\leq n$) corresponds to the $j$-th column of
$v^{(-\infty,q]}_{\nu}$ from the right. By Lemma \ref{babb commutes}
(and Remark \ref{remark on babb commutes}), we have
{\allowdisplaybreaks
\begin{align*}
(A_1,\ldots,A_{m-1},A_m,A_{m+1},\ldots,A_{m+n})& \equiv (A_1,\ldots,A_{m-1},A_{m+1},A_{m},\ldots,A_{m+n}) \\
& \equiv (A_1,\ldots,A_{m+1},A_{m-1},A_{m},\ldots,A_{m+n})\\
& \ \ \ \ \ \ \ \ \ \ \ \ \ \ \ \ \ \ \ \ \ \ \  \vdots \\
& \equiv (A_{m+1},A_1\ldots,A_m,A_{m+2},\ldots,A_{m+n}).
\end{align*}}
Note that  we can identify each matrix given above with an element in a $\gl_\infty$-crystal (a mixed tensor product of $\E$ and $\E^\vee$'s), and hence consider $\gl_\infty$-crystal equivalence between them.

Repeating the above process, we conclude that
$$(A_1,\ldots,A_m,A_{m+1},\ldots,A_{m+n})\equiv
(A_{m+1},\ldots,A_{m+n},A_1,\ldots,A_m).$$

Since the righthand side
of the above equivalence is the image of
$v^{(-\infty,q]}_{\nu}\otimes u^{[p,\infty)}_{\mu}$ in $\left( \E^n\right)^\vee\otimes\E^m \subset \bM_{[m+n],\Z}$ under
$\tau_\nu\times\sigma_\mu$, we have
$$u^{[p,\infty)}_{\mu}\otimes v^{(-\infty,q]}_{\nu}\equiv v^{(-\infty,q]}_{\nu}\otimes u^{[p,\infty)}_{\mu}.$$
By Proposition \ref{connectedness of bmn}, it follows that
$\cB_{\mu}\otimes\cB_{\nu}^\vee\simeq \cB_{\nu}^\vee\otimes
\cB_{\mu}$. \qed\vskip 3mm

From now on, we write $\cB_{\mu,\nu}=\cB_{\mu}\otimes\cB_\nu^\vee$ for $\mu,\nu\in\cP$.

\begin{prop}\label{equi of bmn} For $\mu,\nu,\sigma,\tau\in\cP$,
$\cB_{\mu,\nu} \simeq \cB_{\sigma,\tau}$ if and only if
$(\mu,\nu)=(\sigma,\tau)$.
\end{prop}
\pf Suppose that $S\otimes T\equiv S'\otimes T'$ for $S\otimes T\in
\cB_{\mu}\otimes\cB_{\nu}^\vee$ and $S'\otimes T'\in
\cB_{\sigma}\otimes\cB_{\tau}^\vee$. By Remark \ref{remark on bmn},
we may assume that there exist $s\ll 0$ and $t\gg 0$ such that the
entries in $S$ and $S'$ are less than $s$, and the entries in $T$
and $T'$ are less than $t^\vee$. Then it follows that $S$ is
$\gl_{(-\infty,s]}$-equivalent to $S'$, which implies that $S=S'$
and $\mu=\sigma$. Similarly, we have $T=T'$ and $\nu=\tau$.
\qed\vskip 3mm

We define ${\rm wt}_{[n]}(A)= \sum_{1\leq i\leq
n}\left(\sum_{j\in\Z}a_{ij}\right)\epsilon_i$ for $A\in\E^n$.  Then $\E^n$ ($n\geq 2$) is
column admissible, and it is a regular
$\gl_{[n]}$-crystal with respect to $\tE_k, \tF_k$ ($1\leq k\leq n-1$) and ${\rm
wt}_{[n]}$. By Lemma \ref{commutativity}, $\E^n$ is a
$(\gl_\infty,\gl_{[n]})$-bicrystal, that is, the operators
$\te_i,\tf_i$ $(i\in\Z)$ on $\E^{n}\cup\{{\bf 0}\}$ commute with
$\tE_j, \tF_j$  $(1\leq j\leq n-1)$.

For $k\in [n]$ and $\lambda\in \Z_+^n$, we put 
\begin{equation}
\begin{split}
\omega_{k}&=\epsilon_{1}+\cdots+\epsilon_{k}, \\
\omega_{\lambda}&=\lambda_1\epsilon_{1}+\cdots+\lambda_{n}\epsilon_{n}.
\end{split}
\end{equation}
We denote by $\B_{[n]}(\pm \omega_{\lambda})$ the crystal base
of the irreducible $U_q(\gl_{[n]})$-module with highest (resp. lowest) weight
vector $u^{[n]}_{\pm\omega_{\lambda}}$ of weight $\pm\omega_{\lambda}$.
When $\lambda\in\cP$, $\B_{[n]}(\omega_{\lambda})$ can be realized as $SST_{[n]}(\lambda)$ \cite{KN}.
Here we assume that $\B_{[1]}(\pm\omega_\lambda)=\{\,u^{[1]}_{\pm\omega_{\lambda}}\,\}$ with ${\rm wt}_{[1]}(u^{[1]}_{\pm\omega_{\lambda}})=\pm\lambda\epsilon_1$ for $\lambda\in \Z^1_+=\Z$.

Note that for $A\in\E^n$ ($n\geq 2$), the $j$-th column $A^j$ of $A$ ($j\in \Z$) is $\gl_{[n]}$-equivalent to  an element in the trivial crystal or $\B_{[n]}(\omega_k)$ for some
$1\leq k\leq n$, and it is non-trivial for only finitely many $j$'s.

\begin{prop}\label{duality on En}
As a $(\gl_\infty,\gl_{[n]})$-bicrystal, we have
$$\E^n \simeq \bigsqcup_{\substack{\mu\in\cP \\ \mu_1\leq n}}\cB_{\mu}\times \B_{[n]}(\omega_{\mu'}).$$
\end{prop}
\pf We may assume that $n\geq 2$. Let $\mu\in\cP$ be given with $\mu_1\leq n$. Consider the image of $u^{[1,\infty)}_{\mu}$ under  $\sigma_\mu$ (\ref{sigma'}), say $A_\mu=(a_{ij})$. Recall that $A$ is of the form;
\begin{equation*}
a_{ij}=
\begin{cases}
1, & \text{if $1 \leq j\leq \mu'_{n-i+1}$}, \\
0, & \text{otherwise}.
\end{cases}
\end{equation*}
Then it is straightforward to see that  $A_\mu$ is $\gl_{[n]}$-equivalent to the lowest weight vector in $\B_{[n]}(\omega_{\mu'})$, that is, $S_{w_{[n]}}u^{[n]}_{\omega_{\mu'}}$, where $w_{[n]}$ is the longest element in the Weyl group of $\gl_{[n]}$. This implies that $C(A_\mu)$ the connected component in $\E^n$ including $A_\mu$ is isomorphic to $\cB_{\mu}\times \B_{[n]}(\omega_{\mu'})$ as a $(\gl_\infty,\gl_{[n]})$-bicrystal. Note that for $\mu,\nu\in \cP$, $C(A_\mu)=C(A_\nu)$   if and only if $\mu=\nu$.

Suppose that $B\in \E^n$ is given. Choose an interval $[p,q]$ in $\Z$ such that all non-zero entries of $B$ are placed in its $[n]\times [p,q]$ submatrix. Since $\bM_{[n],[p,q]}$ is a $(\gl_{[n]},\gl_{[p,q]})$-bicrystal, $B$ is connected to $B'=(b'_{ij})$, which is of the following form;
\begin{equation*}
b'_{ij}=
\begin{cases}
1, & \text{if $p \leq j\leq \mu'_{n-i+1}+p-1$}, \\
0, & \text{otherwise},
\end{cases}
\end{equation*}
for some $\mu\in \cP$ with $\mu_1\leq n$ (see \cite[Theorem
4.5]{K07}). Now, we have $S_{w}B'=A_\mu$ for $w\in W$ such that $w({\rm wt}(B'))={\rm wt}(A_\mu)$, where $S_w$ is defined with respect to $\te_k$ and $\tf_k$ ($k\in\Z$). Hence   $C(B')=C(A_\mu)$. This completes the proof. \qed\vskip 3mm

Considering the dual crystal of $\E^n$, we have
\begin{equation}
\left(\E^{n}\right)^\vee\simeq \bigsqcup_{\substack{\mu\in\cP \\ \mu_1\leq n}}\cB_{\mu}^\vee\times \B_{[n]}(-\omega_{\mu'}).
\end{equation}

\subsection{Highest weight crystals}
Let $\F$ be the subset of $\bM_{\{1\},\Z}$ consisting of
$(a_{1j})_{j\in\mathbb{Z}}$ such that $a_{1j}=1$ for
$j\ll 0$ and $a_{1j}=0$ for $j\gg 0$.
For
$A=(a_{1j})_{j\in\mathbb{Z}}\in\F$, define
\begin{equation}\label{weight}
{\rm wt}(A)=\Lambda_0+\sum_{j>0}a_{1j}\epsilon_j+\sum_{j\leq
0}(a_{1j}-1)\epsilon_j.
\end{equation}
Then $\F$ is a regular $\gl_{\infty}$-crystal with respect to $\te_k, \tf_k$ ($k\in\Z$) and ${\rm wt}$, and
there exists a strict embedding
\begin{equation}\label{iota_i}
\iota_i : \B(\Lambda_i) \longrightarrow \F
\end{equation}
for $i\in\Z$, where  the highest weight vector $u_{\Lambda_i}$ of $\B(\Lambda_i)$ is mapped to the unique element
of weight $\Lambda_i$, that is, $\iota_i(u_{\Lambda_i})=(a_{1j})_{j\in\Z}$ with $a_{1j}=1$ for $j\leq i$ and $a_{1j}=0$ otherwise.
Then, we have
\begin{equation}\label{F1 decomposition}
\F\simeq \bigsqcup_{i\in\mathbb{Z}}\B(\Lambda_i).
\end{equation}
Recall that $\F$ is the crystal base
of the $q$-deformed Fock space representation, which can be realized as the space
of semi-infinite wedge vectors \cite{MM,U}.

For $n\geq 1$, let $\F^n$ be the set of matrices $A=(a_{ij})$ in $\bM_{[n],\Z}$ such that $A_i\in \F\subset \bM_{\{i\},\Z}$ for $i\in[n]$. Then $\F^n$ is row admissible and can be identified with $\F^{\otimes n}$ as  a $\gl_\infty$-crystal, where $A$ is identified with $A_1\otimes\cdots\otimes A_n$. Also for $\lambda\in\Z^n_+$, we may regard $\B(\Lambda_\lambda)\subset \B(\Lambda_{\lambda_n})\otimes\cdots\otimes\B(\Lambda_{\lambda_1})$ by identifying $u_{\Lambda_\lambda}$ with $u_{\Lambda_{\lambda_n}}\otimes\cdots\otimes u_{\Lambda_{\lambda_1}}$. Composing with (\ref{iota_i}), we have a strict embedding
\begin{equation}\label{embedding of HW}
\iota_\lambda=\iota_{\lambda_n}\times\cdots\times\iota_{\lambda_1} : \B(\Lambda_\lambda) \longrightarrow   \F^n.
\end{equation}
Taking dual crystals in (\ref{iota_i}) and (\ref{embedding of HW}), we have embeddings $\iota^\vee_i$ and $\iota^\vee_\lambda$, respectively.

On the other hand, define
\begin{equation}
{\rm wt}_{[n]}(A)=\sum_{1\leq i\leq n}\sum_{j>0}a_{ij}\epsilon_i + \sum_{1\leq i\leq n}\sum_{j\leq 0}(a_{ij}-1)\epsilon_i.
\end{equation}
for  $A\in\F^n$.
Then $\F^n$ ($n\geq 2$) is column admissible and it is a regular $\gl_{[n]}$-crystal  with respect to $\tE_k, \tF_k$ and ${\rm wt}_{[n]}$. Hence, $\F^n$ is a $(\gl_\infty,\gl_{[n]})$-bicrystal  by Lemma \ref{commutativity}.

Note that for $A\in\F^n$ ($n\geq 2$), the $j$-th column $A^j$ of $A$ ($j\in \Z$) is $\gl_{[n]}$-equivalent to  an element in the trivial crystal or $\B_{[n]}(\pm\omega_k)$ for some
$1\leq k\leq n$, and it is non-trivial for only finitely many $j$'s.

The following theorem is a crystal version of the
$(\gl_\infty,\gl_{[n]})$-duality on the level $n$ fermionic Fock space
\cite{Fr}.
\begin{prop}[cf.\cite{K09}]\label{duality-1}
As
a $(\gl_{\infty},\gl_{[n]})$-bicrystal, we have
$$\F^{n}\simeq\bigsqcup_{\lambda\in\Z_+^n}\B(\Lambda_{\lambda})\times \B_{[n]}(\omega_{\lambda}).$$
\end{prop}
If we consider $\left(\F^{n}\right)^\vee$, then we obtain
$$\left(\F^{n}\right)^\vee\simeq\bigsqcup_{\lambda\in\Z_+^n}\B(-\Lambda_{\lambda})\times \B_{[n]}(-\omega_{\lambda}).$$
As in the case of $\left(\E^n\right)^\vee$, we may view $\left(\F^{n}\right)^\vee=\{\,A^\vee\,|\,A\in\F^n\,\}$.

\subsection{Littlewood-Richardson coefficients}
Let $x=\{\,x_1,x_2,x_3,\ldots\,\}$ be the set of formal commuting
variables. Let $Sym$ be the ring of symmetric functions in $x$.  For
$k\geq 1$, denote  by $e_k(x)$, $h_k(x)$ and $p_k(x)$ the $k$-th
elementary, complete and power sum symmetric functions in $x$,
respectively. It is well-known that $\{\,e_k(x)\,|\,k\geq 1\,\}$ and
$\{\,h_k(x)\,|\,k\geq 1\,\}$ are algebraically independent over
$\mathbb{Z}$ in $Sym$ and $\{\,p_k(x)\,|\,k\geq 1\,\}$ is
algebraically independent over $\mathbb{Q}$ in
$Sym_{\mathbb{Q}}=\mathbb{Q}\otimes_{\Z}Sym$.

For $\lambda\in\cP$, let $s_{\lambda}(x)$ be the Schur function in $x$
corresponding to $\lambda$ \cite{Mac95}. The Littlewood-Richardson
coefficients $c^{\lambda}_{\mu\,\nu}$ for $\lambda,\mu,\nu\in\cP$ are
defined by
\begin{equation}
s_{\mu}(x)s_{\nu}(x)=\sum_{\lambda}c^{\lambda}_{\mu\,\nu}s_{\lambda}(x).
\end{equation}
For $n\geq 1$, let $x_{[n]}=\{\,x_1,\ldots,x_n\,\}$. For
$\lambda\in\cP$ with $\ell(\lambda)\leq n$, let
$s_{\lambda}(x_{[n]})$ be the corresponding Schur polynomial in
$x_{[n]}$. Put ${\rm
ch}\B_{[n]}(\omega_\lambda)=\sum_{T\in \B_{[n]}(\omega_\lambda)}x^{T}_{[n]}$, where
$x_{[n]}^T=\prod_i x_i^{m_i}$ for $T\in \B_{[n]}(\omega_\lambda)$
with ${\rm wt}_{[n]}(T)=\sum_im_i\epsilon_i$. Then we have ${\rm
ch}\B_{[n]}(\omega_\lambda)=s_{\lambda}(x_{[n]})$.

Hereafter, for a crystal $B$ and a non-negative integer $m$, we
denote by $B^{\oplus m}$ the disjoint union $B_1\sqcup\cdots\sqcup
B_m$ with $B_i\simeq B$, where $B^{\oplus 0}$ means the empty set.
\begin{prop}\label{LR for BmBn}
For $\mu,\nu\in \cP$, we have
$$\cB_{\mu}\otimes\cB_{\nu}\simeq \bigsqcup_{\lambda\in\cP}\cB_\lambda^{\oplus c^\lambda_{\mu \nu}}.$$
\end{prop}
\pf Let $m,n$ be positive integers. Put $[n]+m=\{\,m+1,\ldots,m+n\,\}$.
Then $\gl_{[m]}\oplus\gl_{[n]+m}$ is a subalgebra of $\gl_{[m+n]}$.
By Proposition \ref{duality on En}, we have
\begin{equation}\label{Em times En}
\E^m\otimes\E^n\simeq \bigsqcup_{\substack{\mu,\nu\in\cP \\
\mu_1\leq m,\,\nu_1\leq n}} \left(\cB_\mu\otimes\cB_\nu \right)\times \left(\B_{[m]}(\omega_{\mu'})\times \B_{[n]+m}(\omega_{\nu'}) \right),
\end{equation}
as a $(\gl_\infty,\gl_{[m]}\oplus\gl_{[m]+n})$-bicrystal. On the other hand, we have
\begin{equation}\label{Em+n}
\E^{m+n}\simeq \bigsqcup_{\substack{\lambda\in\cP \\ \lambda_1\leq m+n}} \cB_\lambda\times \B_{[m+n]}(\omega_{\lambda'}),
\end{equation}
as a $(\gl_\infty,\gl_{[m+n]})$-bicrystal. Since
$s_{\lambda'}(x_{[m+n]})=\sum_{\mu'\nu'}c^{\lambda'}_{\mu'
\nu'}s_{\mu'}(x_{[m]})s_{\nu'}(x_{[n]+m})$ and $c^{\lambda'}_{\mu'
\nu'}=c^{\lambda}_{\mu \nu}$, we have as a
$\gl_{[m]}\oplus\gl_{[m]+n}$-crystal,
$$\B_{[m+n]}(\omega_{\lambda'})\simeq
\bigsqcup_{\mu',\nu'}\B_{[m]}(\omega_{\mu'})\times
\B_{[n]+m}(\omega_{\nu'})^{\oplus {c^\lambda_{\mu \nu}}}.$$ Since
$\E^{m+n}\simeq \E^m\otimes\E^n$ as a
$(\gl_\infty,\gl_{[m]}\oplus\gl_{[m]+n})$-bicrystal, we obtain the
required decomposition of $\cB_{\mu}\otimes\cB_{\nu}$ by comparing
(\ref{Em times En}) and (\ref{Em+n}). \qed\vskip 3mm

Let $m,n\geq 1$ be given. Recall that for $\mu\in\Z_+^m$, ${\rm ch}\B_{[m]}(\omega_\mu)=s_{\mu}(x_{[m]})=(x_1\cdots x_m)^{-p} s_{\mu+(p^m)}(x_{[m]})$ is the Laurent Schur polynomial corresponding to $\mu$, where $p$ is a non-negative integer such that $\mu+(p^m)\in\cP$.

For $\lambda\in\Z_+^{m+n}$, $\mu\in\Z_+^m$ and $\nu\in\Z_+^n$, we define
\begin{equation}
c^\lambda_{\mu \nu}=c^{\lambda+(p^{m+n})}_{\mu+(p^m) \nu+(p^n)},
\end{equation}
where $p$ is a non-negative integer such that $\lambda+(p^{m+n}),\mu+(p^m),
\nu+(p^n)\in\cP$. Note that
$s_{\lambda}(x_{[m+n]})=\sum_{\mu,\nu}c^\lambda_{\mu
\nu}s_{\mu}(x_{[m]})s_{\nu}(x_{[n]+m})$ and $c^\lambda_{\mu \nu}$
does not depend on $p$.

\begin{prop}\label{LR HW}
For $\mu\in\Z_+^m$ and $\nu\in\Z_+^n$,
$$\B(\Lambda_\mu)\otimes \B(\Lambda_\nu)\simeq \bigsqcup_{\lambda\in\Z_+^{m+n}}\B(\Lambda_{\lambda})^{\oplus c^\lambda_{\mu \nu}}.$$
\end{prop}
\pf  The proof is almost the same as that of Proposition \ref{LR for BmBn}. Here we compare the $(\gl_{\infty},\gl_{[m]}\oplus\gl_{[n]+m})$-bicrystal decompositions of $\F^m\otimes\F^n$ and $\F^{m+n}$. \qed\vskip 3mm

\begin{rem}{\rm
Note that there are infinitely many connected components in $\B(\Lambda_\mu)\otimes\B(\Lambda_\nu)$, but the multiplicity of each connected component is finite.
}
\end{rem}

\subsection{Tensor product}

Let us end this section with introducing another family of regular connected $\gl_\infty$-crystals.

\begin{prop}\label{Blmn is connected} For $\mu,\nu\in\cP$ and $\Lambda\in P^+$,
$\cB_{\mu,\nu}\otimes \B(\Lambda)$ is connected.
\end{prop}
\pf Suppose that $S\otimes T\otimes U \in \cB_{\mu,\nu}\otimes
\B(\Lambda)$ is given. By Remark \ref{remark on bmn}, there exist
$i_1,\ldots,i_r\in\Z$  such that $\te_{i_1}\cdots\te_{i_r}(S\otimes
T)=u^{[p,\infty)}_{\mu}\otimes v^{(-\infty,q]}_{\nu}$ for some
$p<q$. By tensor product rule of crystals, we have
$$\te_{i_1}^{m_1}\cdots\te_{i_r}^{m_r}(S\otimes T\otimes U)=u^{[p,\infty)}_{\mu}\otimes v^{(-\infty,q]}_{\nu}\otimes U'$$
for some $m_1,\ldots,m_r\geq 1$ and $U'\in\B(\Lambda)$.  We may
assume that $p\ll 0$ and $q\gg 0$ so that
$\te_{j_1}\cdots\te_{j_s}U'=u_{\Lambda}$ for some $j_1,\ldots,j_s\in
[p+\ell(\mu)+1,q-\ell(\nu)-1]$. Since
$\widetilde{x}_{j_t}\left(u^{[p,\infty)}_{\mu}\otimes
v^{(-\infty,q]}_{\nu} \right)={\bf 0}$ for $1\leq t\leq s$ and
$x=e,f$, we get
$$\te_{j_1}\cdots\te_{j_s}\left(u^{[p,\infty)}_{\mu}\otimes v^{(-\infty,q]}_{\nu}\otimes U'\right)
=u^{[p,\infty)}_{\mu}\otimes v^{(-\infty,q]}_{\nu}\otimes u_{\Lambda}.$$
Since $p$ (resp. $q$) can be arbitrarily small (resp. large), we
conclude that $\cB_{\mu,\nu}\otimes \B(\Lambda)$ is connected.
\qed\vskip 3mm

\begin{prop}\label{equi of bmnBLambda} For $\mu,\nu,\sigma,\tau\in\cP$ and $\Lambda,\Lambda'\in P^+$,
$\cB_{\mu, \nu}\otimes \B(\Lambda)\simeq \cB_{\sigma, \tau}\otimes
\B(\Lambda')$ if and only if
$(\mu,\nu,\Lambda)=(\sigma,\tau,\Lambda')$.
\end{prop}
\pf Suppose that $\cB_{\mu, \nu}\otimes \B(\Lambda)\simeq
\cB_{\sigma, \tau}\otimes \B(\Lambda')$. Let $S\otimes T\otimes U\in
\cB_{\mu, \nu}\otimes \B(\Lambda)$ be equivalent to $S'\otimes
T'\otimes U'\in \cB_{\sigma, \tau}\otimes \B(\Lambda')$.  Applying
suitable $\te_k$'s, we assume that there exist $s\ll 0$ and $t\gg 0$
such that the entries in $S$ and $S'$ are less than $s$, and the
entries in $T$ and $T'$ are less than $t^\vee$ (see Remark
\ref{remark on bmn} and the proof of Proposition \ref{Blmn is
connected}). We may further assume that $\widetilde{x}_k
U=\widetilde{x}_k U'={\bf 0}$ for $k\not\in [s,t]$. By similar
argument as in Proposition \ref{equi of bmn}, we have $S=S'$ and
$T=T'$. Hence $U$ is $\gl_{[s,t]}$-equivalent to $U'$, which implies
that $U=U'$. \qed\vskip 3mm

\section{Realization of extremal weight crystals}

\subsection{Extremal weight crystals for $\gl_\infty$}
Let us briefly recall the crystal bases of the
modified quantized enveloping algebra of $\gl_{\infty}$ and an
extremal weight module over $U_q(\gl_{\infty})$ (see
\cite{Kas94',Kas02} for more details). Let $\widetilde{U}_q(\gl_{\infty})=\bigoplus_{\Lambda\in
P}U_q(\gl_{\infty})a_{\Lambda}$ be the modified quantized
enveloping algebra of $\gl_{\infty}$ and let
\begin{equation}
\B(\widetilde{U}_q(\gl_{\infty}))=\bigsqcup_{\Lambda\in
P}\B(U_q(\gl_{\infty})a_{\Lambda})
\end{equation}
be its crystal base.
It is known that $\B(\widetilde{U}_q(\gl_{\infty}))$ is regular, and
\begin{equation}
\B(U_q(\gl_{\infty})a_{\Lambda})\simeq \B
(\infty)\otimes T_{\Lambda}\otimes\B (-\infty)
\end{equation}
for $\Lambda\in
P$, where $\B(\infty)$ (resp. $\B (-\infty)$) is the crystal
base of the negative (resp. positive) part of $U_q(\gl_{\infty})$, and
$T_{\Lambda}=\{\,t_{\Lambda}\,\}$ is a $\gl_{\infty}$-crystal with
${\rm wt}(t_{\Lambda})=\Lambda$, $\te_i t_{\Lambda}=\tf_{i}
t_{\Lambda}={\bf 0}$, and
$\varepsilon_{i}(t_{\Lambda})=\varphi_{i}(t_{\Lambda})=-\infty$ for
$i\in\Z$.

 An element $b$ of a regular $\gl_{\infty}$-crystal $B$ with
${\rm wt}(b)=\Lambda$ is called {\it extremal} if $\{\,S_wb\,|\,w\in
W \,\}$ satisfies the following conditions; (1) $\te_iS_w b={\bf 0}$
if $\langle w(\Lambda),h_i \rangle\geq 0$, (2) $\tf_iS_w(b)={\bf 0}$
if $\langle w(\Lambda),h_i \rangle\leq 0$.

For $\Lambda\in P $, let
\begin{equation}
\B (\Lambda)=\{\,b\in \B(U_q(\gl_{\infty})a_{\Lambda})\,|\,\text{$b^*$
is extremal}\,\},
\end{equation}
where $\ast$ is the star operation on
$\B(\widetilde{U}_q(\gl_{\infty}))$. Then $\B (\Lambda)$ is the
crystal base of the $U_q(\gl_{\infty})$-module generated by an
extremal weight vector $u_{\Lambda}$ of weight $\Lambda$, which is
called an extremal weight module. Note that (1) $\B(\Lambda)\simeq
\B(w\Lambda)$ for $w\in W $, and (2) if $\Lambda\in P^+$, then
$\B(\Lambda)\simeq \B (\Lambda_{\lambda})$ for some
$\lambda\in\Z_+^n$. From now on, we call $\B(\Lambda)$ simply an
{\it extremal weight crystal}.

\begin{prop}\label{connectedness}
For $\Lambda\in P $, $\B (\Lambda)$ is connected.
\end{prop}
\pf We regard $\B (\Lambda)$ as a subcrystal of $\B (\infty)\otimes
T_{\Lambda}\otimes\B (-\infty)$ and identify
$u_{\Lambda}\in\B(\Lambda)$ with $u_{\infty}\otimes
t_{\Lambda}\otimes u_{-\infty}$, where $u_{\pm \infty}$ is the
highest (resp. lowest) weight vector in $\B(\pm\infty)$. Let $b\in
\B (\Lambda)$ be given. We may assume that $b$ is extremal since any
element in $\B(\Lambda)$ is connected to an extremal one. By the
same argument as in \cite[Theorem 5.1]{Kas02}, $b$ is connected to
$b_1\otimes t_{\Lambda}\otimes u_{-\infty}$, where $\langle {\rm wt}
(b_1),h_i\rangle\geq 0$ for all $i\in\Z$. Since ${\rm wt}
(b_1)=\sum_{i\in\Z}m_i\alpha_i=\sum_{i\in\Z}m_i(\epsilon_i-\epsilon_{i+1})$
with $m_i\in\mathbb{Z}_{\leq 0}$ and $P^+_0=\{0\}$, we have $m_i=0$ for all
$i\in\Z$ and $b_1=u_{\infty}$. Therefore, $\B(\Lambda)$ is
connected. \qed\vskip 3mm

\begin{cor}\label{embedding of BLambda}
For $\Lambda\in P$, $\B(\Lambda)$ is isomorphic to  the connected
component in $\B(\Lambda_+)\otimes\B(-\Lambda_-)$ including
$u_{\Lambda_+}\otimes u_{-\Lambda_-}$.
\end{cor}
\pf Recall that there is a strict embedding of regular crystals
$$\B(\Lambda_+)\otimes\B(-\Lambda_-)\longrightarrow \B
(\infty)\otimes T_{\Lambda}\otimes\B (-\infty)$$ sending
$u_{\Lambda_+}\otimes u_{-\Lambda_-}$ to $u_{\infty}\otimes
t_{\Lambda}\otimes u_{-\infty}$. Since $\B(\Lambda)\simeq
C(u_{\infty}\otimes t_{\Lambda}\otimes u_{-\infty})$ by Proposition
\ref{connectedness}, we have $\B(\Lambda)\simeq
C(u_{\Lambda_+}\otimes u_{-\Lambda_-})\subset
\B(\Lambda_+)\otimes\B(-\Lambda_-)$. \qed\vskip 3mm

\subsection{Realization of extremal weight crystals}
\begin{lem}\label{level 1-1 decomposition}
For $i,j\in\Z$ $(i\leq j)$, we have  {\allowdisplaybreaks
\begin{equation*}
\begin{split}
&\B(\Lambda_i)\otimes\B(-\Lambda_j)\simeq \bigsqcup_{a\geq 0}\cB_{(1^{a}),(1^{a+j-i})}, \\
&\B(\Lambda_j)\otimes\B(-\Lambda_i)\simeq \bigsqcup_{a\geq 0}\cB_{(1^{a+j-i}),(1^{a})}. \\
\end{split}
\end{equation*}}
\end{lem}
\pf Let us prove the first isomorphism. The second one is obtained
by considering the dual crystals on both sides of the first
isomorphism.

Suppose that $S\otimes T\in \B(\Lambda_i)\otimes\B(-\Lambda_j)$ is
given.  Applying $\te_k$'s, we assume that
$S=u_{\Lambda_i}$. Then applying $\tf_k$'s ($k\neq i$),
$u_{\Lambda_i}\otimes T$ is connected to $u_{\Lambda_i}\otimes T'$
such that $\tf_k\left(u_{\Lambda_i}\otimes T'\right)={\bf 0}$ for
all $k\neq i$.

Let $A=(a_{lk})$ be the image of $u_{\Lambda_i}\otimes T'$ in $\bM_{[2],\Z}$ under the strict embedding $\iota_i\times \iota_j^\vee : \B(\Lambda_i)\otimes \B(-\Lambda_j)\rightarrow \F\otimes \F^\vee$. Then there exists $a\geq
0$ such that
\begin{equation*}\label{ext vector}
\begin{split}
\begin{cases}
a_{1k}= 1 & \text{if and only if $k\leq i$},\\
a_{2k}= 1 & \text{if and only if $i-a+1\leq k\leq i$ or $j+a+1\leq k$}.\\
\end{cases}
\end{split}
\end{equation*}
Let $B=(b_{lk})\in\E\otimes \E^\vee\subset \bM_{[2],\Z}$ be such that
\begin{equation*}
\begin{cases}
b_{1k}= 1 & \text{if and only if $i-a+1\leq k\leq i$},\\
b_{2k}= 0 & \text{if and only if $i+1\leq k\leq j+a$}.
\end{cases}
\end{equation*}
Note that $C(B)\simeq \cB_{(1^{a}),(1^{a+j-i})}$. For $p<q\in\Z$, let
$\pi_{[p,q]} : \bM_{[2],\Z}\rightarrow \bM_{[2],[p,q]}$ be the map
sending a matrix to its $[2]\times [p,q]$ submatrix. Assume that $p\ll 0
\ll q$. Then $A$ is column admissible and $\tF_1^{\rm max} \pi_{[p,q]}(A)=\pi_{[p,q]}(B)$.
By Lemma
\ref{admissibility}, we have
\begin{equation*}
\widetilde{x}_{i_1}\cdots\widetilde{x}_{i_r} \pi_{[p,q]}(A) \neq {\bf 0} \Longleftrightarrow \widetilde{x}_{i_1}\cdots\widetilde{x}_{i_r} \pi_{[p,q]}(B)\neq {\bf 0},
\end{equation*}
and hence
\begin{equation*}
\widetilde{x}_{i_1}\cdots\widetilde{x}_{i_r} A \neq {\bf 0} \Longleftrightarrow \widetilde{x}_{i_1}\cdots\widetilde{x}_{i_r} B\neq {\bf 0}
\end{equation*}
for $r\geq 1$ and $p\leq i_1,\ldots,i_r\leq q-1$, where $x=e$ or $f$ for each $i_k$.
Since $p$ and $q$ are arbitrary and ${\rm wt}(A)={\rm wt}(B)$, we have $A\equiv B$ or $u_{\Lambda_i}\otimes T'\equiv B$, which implies that $C(u_{\Lambda_i}\otimes T')\simeq \cB_{(1^{a}),(1^{a+j-i})}$.

Conversely, for each $a\geq 0$, there exists a unique $T\in \B(-\Lambda_j)$ such that  
$C(u_{\Lambda_i}\otimes T)\simeq\cB_{(1^{a}),(1^{a+j-i})}$ since the construction of $B$ is reversible.  This
completes the proof. \qed\vskip 3mm

\begin{lem}\label{commuting relation-1}
For $i\in\Z$ and $k\geq 0$, we have  {\allowdisplaybreaks
\begin{equation*}
\begin{split}
&\B(\Lambda_i)\otimes\cB_{(1^k)}\simeq \bigsqcup_{a=0}^k \cB_{(1^a)}\otimes \B(\Lambda_{i+k-a}),\\
&\B(\Lambda_i)\otimes\cB_{(1^k)}^\vee\simeq \bigsqcup_{a=0}^{k} \cB_{(1^a)}^\vee\otimes \B(\Lambda_{i-k+a}).
\end{split}
\end{equation*}}
\end{lem}
\pf The proof is similar to that of Lemma \ref{level 1-1
decomposition}. Let us prove the first isomorphism. Suppose that
$S\otimes T\in \B(\Lambda_i)\otimes\cB_{(1^k)}$ is
given.  Applying $\te_k$'s, we assume that $S=u_{\Lambda_i}$.

Let $A=(a_{lj})$ be the image of $u_{\Lambda_i}\otimes T$ in
$\bM_{[2],\Z}$ under the strict embedding $\iota_i\times \sigma_k : \B(\Lambda_i)\otimes\B_{(1^k)}\rightarrow \F\otimes \E$.  Applying suitable
$\widetilde{x}_s$'s for $s\neq i$ and $x=e,f$, we may assume that
$a_{2j}=1$ if and only if $i-a+1 \leq j\leq i-a+k$ for some $a\geq
0$. Let $B=(b_{lj})\in \E\otimes \F\subset \bM_{[2],\Z}$ be such
that
\begin{equation*}
\begin{cases}
b_{1j}= 1 &   \text{if and only if $i-a+1\leq j\leq i$},\\
b_{2j}= 1 &   \text{if and only if $j\leq i-a+k$}.
\end{cases}
\end{equation*}
Note that  $C(B)\simeq \cB_{(1^a)}\otimes \B(\Lambda_{i-a+k})$. Choose $p\ll 0 \ll q$.
Then $A$ is column admissible and $\tF_1^{\rm max} \pi_{[p,q]}(A)=\pi_{[p,q]}(B)$.
By Lemma
\ref{admissibility}, we have
\begin{equation*}
\widetilde{x}_{i_1}\cdots\widetilde{x}_{i_r} \pi_{[p,q]}(A) \neq {\bf 0} \Longleftrightarrow \widetilde{x}_{i_1}\cdots\widetilde{x}_{i_r} \pi_{[p,q]}(B)\neq {\bf 0},
\end{equation*}
and hence
\begin{equation*}
\widetilde{x}_{i_1}\cdots\widetilde{x}_{i_r} A \neq {\bf 0} \Longleftrightarrow \widetilde{x}_{i_1}\cdots\widetilde{x}_{i_r} B\neq {\bf 0}
\end{equation*}
for $r\geq 1$ and $p\leq i_1,\ldots,i_r\leq q-1$, where $x=e$ or $f$ for each $i_s$.
Since $p$ and $q$ are arbitrary and ${\rm wt}(A)={\rm wt}(B)$, we have $A\equiv B$ or $u_{\Lambda_i}\otimes T\equiv B$, which implies that $C(u_{\Lambda_i}\otimes T)\simeq \cB_{(1^a)}\otimes \B(\Lambda_{i-a+k})$.

Conversely, for $a\geq 0$, we can find a unique
$T\in \B_{(1^k)}$  such that  $C(u_{\Lambda_i}\otimes
T)\simeq \cB_{(1^a)}\otimes \B(\Lambda_{i-a+k})$ since the construction of $B$ is reversible. This establishes
the first isomorphism.

The second isomorphism can be proved by modifying the above argument.
\qed\vskip 3mm

\begin{prop}\label{level m-n decomposition} For $m\geq n\geq 0$,
a connected component in $\F^m\otimes \left(\F^n \right)^\vee$ is
isomorphic to $\cB_{\mu, \nu}\otimes \B(\Lambda)$ for some
$\mu,\nu\in\cP$ and $\Lambda\in P^+_{m-n}$.
\end{prop}
\pf We claim that each $A\in \F^m\otimes \left(\F^n \right)^\vee$ is
equivalent to an element in
$\E^n\otimes\left(\E^n\right)^\vee\otimes\F^{m-n}$. Then it follows  from
Proposition \ref{duality on En}, \ref{duality-1} and \ref{Blmn is
connected} that $C(A)\simeq\cB_{\mu,
\nu}\otimes \B(\Lambda)$ for some $\mu,\nu\in\cP$ and $\Lambda\in
P^+_{m-n}$.

We use induction on $m+n$. Suppose that $m+n=2$. If $m=n=1$,  then
it is clear by Lemma \ref{level 1-1 decomposition}. If $m=2$ and
$n=0$, then it follows from Proposition \ref{duality-1}. Suppose
that $m+n\geq 3$. Let $A=A_1\otimes\cdots\otimes A_{m+n}\in
\F^m\otimes \left(\F^n \right)^\vee$ be given, where $A_i\in\F$ and
$A_{m+j}\in\F^\vee$ for $1\leq i\leq m$ and $1\leq j\leq n$.
Consider $A_m\otimes A_{m+1}\in  \F\otimes\F^\vee$. By (\ref{F1
decomposition}) and Lemma \ref{level 1-1 decomposition}, $A_m\otimes
A_{m+1}$ is equivalent to some $A'_m\otimes A'_{m+1}\in
\E\otimes \E^\vee$. Applying Lemma \ref{commuting relation-1} to
$A_1\otimes\cdots\otimes A_{m-1}$ and $A'_m\otimes A'_{m+1}$
repeatedly, we can say that $A$ is equivalent to some
$B=B_1\otimes\cdots\otimes B_{m+n}$ in $\E\otimes \E^\vee\otimes
\F^{m-1}\otimes \left(\F^{n-1}\right)^\vee$. By induction
hypothesis, $B_{3}\otimes\cdots\otimes B_{m+n}$ is equivalent to
some $B'_{3}\otimes\cdots\otimes B'_{m+n}$ in
$\E^{n-1}\otimes\left(\E^{n-1}\right)^\vee\otimes\F^{m-n}$.
Finally, by Lemma \ref{babb commutes}, $B_1\otimes B_2\otimes
B'_3\otimes\cdots\otimes B'_{2n}\in \E\otimes \E^\vee\otimes
\E^{n-1}\otimes\left(\E^{n-1}\right)^\vee$ is equivalent to an
element in $\E^{n}\otimes\left(\E^{n}\right)^\vee$. Therefore $A$ is
equivalent to an element in
$\E^n\otimes\left(\E^n\right)^\vee\otimes\F^{m-n}$. This completes
the induction. \qed\vskip 3mm

\begin{thm}\label{realization thm}
For $\Lambda\in P_{\ell}$ $(\ell\in \Z)$, there exist unique
$\mu,\nu\in\cP$  and $\Lambda'\in P^+_{|\ell|}$ such that
$$
\B(\Lambda)\simeq
\begin{cases}
\cB_{\mu, \nu}\otimes \B(\Lambda'), & \text{if $\ell\geq 0$}, \\
\B(-\Lambda')\otimes \cB_{\mu, \nu}, & \text{if $\ell\leq 0$}.
\end{cases}
$$
\end{thm}
\pf The first isomorphism follows immediately  from Proposition \ref{duality-1}, Corollary
\ref{embedding of BLambda} and Proposition \ref{level m-n
decomposition}. If $\ell\leq 0$, then $\B(\Lambda)$ is embedded into
$\F^m\otimes \left(\F^n \right)^\vee$ for some $m,n$ with
$m-n=\ell$. Since $\B(\Lambda)^\vee$ is embedded into $\F^n\otimes
\left(\F^m \right)^\vee$, it is isomorphic to $\cB_{\mu, \nu}\otimes
\B(\Lambda')$ for some $\Lambda'\in P^+_{|\ell|}$ and
$\mu,\nu\in\cP$. Hence, $\B(\Lambda)\simeq \B(-\Lambda')\otimes
\cB_{\nu, \mu}$. The uniqueness follows from Proposition \ref{equi
of bmnBLambda}. \qed\vskip 3mm

\begin{cor}\label{level 0 extremal}
For $\Lambda\in P_0$, there exist unique $\mu,\nu\in\cP$ such that $\B(\Lambda)\simeq \cB_{\mu,\nu}$.
\end{cor}

\begin{rem}{\rm
Combining with a tableaux description of $\B(\Lambda)$ with
$\Lambda\in P^+$ (see for example, \cite{K09'}), we obtain a combinatorial
realization of an extremal weight crystal.}
\end{rem}

\subsection{Pieri rules of extremal weight crystals}
We have the following generalization of Lemma \ref{commuting relation-1}.
\begin{prop}\label{Pieri} For $\lambda\in\Z^n_+$ and $k\geq 1$, we have
\begin{equation*}
\begin{split}
\B(\Lambda_\lambda)\otimes \cB_{(1^k)}&\simeq \bigsqcup_{a=0}^k\bigsqcup_{\substack{\mu\in \Z_+^n \\ (\mu-(\lambda_n^n))/(\lambda-(\lambda_n^n)) : \\ \text{a horizontal strip of length $k-a$}}}\cB_{(1^{a})}\otimes\B(\Lambda_\mu), \\
\B(\Lambda_\lambda)\otimes \cB_{(1^k)}^\vee &\simeq \bigsqcup_{a=0}^k\bigsqcup_{\substack{\nu\in \Z_+^n \\ (\lambda-(\nu_n^n))/(\nu-(\nu_n^n)) : \\ \text{a horizontal strip of length $k-a$}}}\cB_{(1^{a})}^\vee \otimes\B(\Lambda_\nu). \\
\end{split}
\end{equation*}
\end{prop}
\pf First, consider $\B(\Lambda_\lambda)\otimes \cB_{(1^k)}$. Given
$S\otimes T\in \B(\Lambda_\lambda)\otimes \cB_{(1^k)}$, let
$A=(a_{ij})$ be the image of $S\otimes T$ in $\bM_{[n+1],\Z}$ under
the strict embedding $\iota_\lambda\times \sigma_k :
\B(\Lambda_\lambda)\otimes \cB_{(1^k)}\rightarrow \F^n\otimes
\E$. By applying suitable $\widetilde{x}_s$'s
for $x=e,f$ and  $s\in \Z$, we may assume that $A$ is of the
following form;
\begin{itemize}
\item[(1)] for $i\in [n]$, $a_{ij}=1$ if and only if $j\leq \lambda_{n-i+1}$, that is, $S=u_{\Lambda_{\lambda}}$,

\item[(2)] $a_{n+1\, j}=0$ for $j\leq \lambda_n-k$,

\item[(3)] $\sum_{i\in [n+1]}a_{ij}\geq \sum_{i\in
[n+1]}a_{i\,j+1}$ for $j\in [\lambda_n-k+1,\infty)$.
\end{itemize}
Let $a=\sum_{j=\lambda_n-k+1}^{\lambda_n} a_{n+1\, j}$ and let $\mu\in\Z_+^n$ be given by
\begin{equation*}\label{highest weight}
\mu_i=
\begin{cases}
\lambda_1 + \sum_{j>\lambda_1}a_{n+1\,j}, & \text{if $i=1$}, \\
\lambda_i+\sum_{j=\lambda_i+1}^{\lambda_{i-1}}a_{n+1\,j}, & \text{if $2\leq i\leq n$ and $\lambda_i<\lambda_{i-1}$},\\
\lambda_i, & \text{if $2\leq i\leq n$ and $\lambda_i=\lambda_{i-1}$}.
\end{cases}
\end{equation*}
Note that $\mu$ is a well-defined generalized partition by (1) and
(3), and $(\mu-(\lambda_n^n))/(\lambda-(\lambda_n^n))$ is a
horizontal strip of length $k-a$. We denote the matrix of the above form by
$A_{a,\mu}$. Let $B_{a,\mu}=(b_{ij})$ be the image of
$u^{[\lambda_n-k+1,\infty)}_{(1^a)}\otimes u_{\Lambda_\mu}$ in $\bM_{[n+1],\Z}$ under the strict embedding
$\sigma_a\times \iota_\mu : \B_{(1^a)}\otimes \B(\Lambda_\mu)\rightarrow \E\otimes \F^n$. Choose $p\ll \lambda_n-k+1\ll q$. Then $A_{a,\mu}$ is column admissible and
\begin{equation*}
\tF_1^{\rm max}\cdots \tF_{n}^{\rm
max}\pi_{[p,q]}(A_{a,\mu})=\pi_{[p,q]}(B_{a,\mu}),
\end{equation*}
where  $\pi_{[p,q]} : \bM_{[n+1],\Z}\rightarrow \bM_{[n+1],[p,q]}$
is the map sending a matrix to its $[n+1]\times [p,q]$ submatrix.
By the same arguments as in Lemmas \ref{level 1-1 decomposition} and \ref{commuting relation-1}, we have $A_{a,\mu}\equiv B_{a,\mu}$, which implies that $C(A_{a,\mu})\simeq C(B_{a,\mu})\simeq \B_{(1^a)}\otimes
\B(\Lambda_\mu)$.  Conversely, suppose that $0\leq a\leq k$ and
$\mu\in\Z_+^n$ are given where
$(\mu-(\lambda_n^n))/(\lambda-(\lambda_n^n))$ is a horizontal strip
of length $k-a$. Then we can check that there exists a unique $T\in \B_{(1^k)}$ such that $C(u_{\Lambda_\lambda}\otimes T) \simeq \B_{(1^a)}\otimes \B(\Lambda_\mu)$ since the construction of $B_{a,\mu}$ is reversible. This proves
the first isomorphism.

Next, consider $\B(\Lambda_\lambda)\otimes \cB_{(1^k)}^\vee$. Let
$\iota^*_\lambda : \B(\Lambda_\lambda) \rightarrow \F^n$  be the embedding
which sends $u_{\Lambda_\lambda}$ to
$u_{\Lambda_{\lambda_1}}\otimes\cdots\otimes
u_{\Lambda_{\lambda_n}}$ (cf. (\ref{embedding of HW})). Given
$S\otimes T\in \B(\Lambda_\lambda)\otimes \cB_{(1^k)}^\vee$, let
$A=(a_{ij})$ be the image of $S\otimes T$ in $\bM_{[n+1],\Z}$ under
the strict  embedding $\iota^*_\lambda\times \tau_k :
\B(\Lambda_\lambda)\otimes \cB_{(1^k)}^\vee\rightarrow \F^n\otimes
\E^\vee$. By applying suitable
$\widetilde{x}_s$'s for $x=e,f$ and $s\in \Z$, we may assume that
$A$ is of the following form;
\begin{itemize}
\item[($1^*$)] for $i\in [n]$, $a_{ij}=1$ if and only if $j\leq \lambda_{i}$, that is, $S=u_{\Lambda_{\lambda}}$,

\item[($2^*$)] $a_{n+1\,j}=1$ for $j\geq \lambda_1+k+1$,

\item[($3^*$)] $\sum_{i\in [n+1]}a_{ij}\leq \sum_{i\in
[n+1]}a_{ij-1}$ for $j\in (-\infty,\lambda_1+k]$.
\end{itemize}
Let $a=\sum_{j=\lambda_1+1}^{\lambda_1+k}(1- a_{n+1\, j})$ and let $\nu\in\Z_+^n$ be given by
\begin{equation*}\label{highest weight-2}
\nu_i=
\begin{cases}
\lambda_i-\sum_{j=\lambda_{i+1}+1}^{\lambda_i}(1-a_{n+1\,j}), & \text{if $1\leq i\leq n-1$ and $\lambda_{i+1}<\lambda_i$},\\
\lambda_i, & \text{if $1\leq i\leq n-1$ and $\lambda_{i+1}=\lambda_i$},\\
\lambda_n - \sum_{j\leq \lambda_n}(1-a_{n+1\,j}), & \text{if $i=n$}.
\end{cases}
\end{equation*}
Note that $\nu$ is a well-defined generalized partition by ($1^*$)
and ($3^*$), and $(\lambda-(\nu_n^n))/(\nu-(\nu_n^n))$ is a
horizontal strip of length $k-a$. We denote the matrix of the above form by
$A^*_{a,\nu}$. Let $B^*_{a,\nu}=(b_{ij})$ be the image of
$v^{(-\infty,\lambda_1+k]}_{(1^a)}\otimes u_{\Lambda_\nu}$ under
$\tau_a\times \iota^*_\nu$. Choose $p\ll \lambda_1+k \ll q$. Then $A^*_{a,\nu}$ is column admissible and 
\begin{equation*}
\tE_1^{\rm max}\cdots \tE_{n}^{\rm
max}\pi_{[p,q]}(A^*_{a,\nu})=\pi_{[p,q]}(B^*_{a,\nu}).
\end{equation*}
As in Lemmas \ref{level 1-1 decomposition} and \ref{commuting relation-1}, we have $A^*_{a,\nu}\equiv B^*_{a,\nu}$, which implies that $C(A^*_{a,\nu})\simeq C(B^*_{a,\nu})\simeq \B^\vee_{(1^a)}\otimes
\B(\Lambda_\nu)$.
Conversely, suppose that $0\leq a\leq k$ and $\nu\in\Z_+^n$ are
given where $(\lambda-(\nu_n^n))/(\nu-(\nu_n^n))$ is a horizontal
strip of length $k-a$. Then there exists a unique $T\in \B^\vee_{(1^k)}$ such that $C(u_{\Lambda_\lambda}\otimes T) \simeq \cB_{(1^{a})}^\vee\otimes\B(\Lambda_\nu)$ since the construction of $B^*_{a,\mu}$ is also reversible.
This proves the second isomorphism. \qed

\begin{cor}\label{Pieri-2}
Let $\lambda\in\Z_+^n$ and $\mu\in \cP$ be given. Then
\begin{itemize}
\item[(1)] $\B(\Lambda_\lambda)\otimes \cB_{\mu}$ is a finite disjoint union of  $\cB_{\nu}\otimes\B(\Lambda_{\eta})$'s for some $\nu\in \cP$ and $\eta\in\Z_+^n$ such that $|\nu|=a\leq |\mu|$ and $(\eta-(\lambda_n^n))/(\lambda-(\lambda_n^n))$ is a skew Young diagram of size $|\mu|-a$,

\item[(2)] $\B(\Lambda_\lambda)\otimes \cB_{\mu}^\vee$ is a finite disjoint union of  $\cB_{\nu}^\vee\otimes\B(\Lambda_{\eta})$'s for some $\nu\in \cP$ and $\eta\in\Z_+^n$ such that $|\nu|=a\leq |\mu|$ and $(\lambda-(\eta_n^n))/(\eta-(\eta_n^n))$ is a skew Young diagram of size $|\mu|-a$.
\end{itemize}
\end{cor}
\pf It follows immediately from Proposition \ref{duality on En} and Proposition \ref{Pieri}.
\qed

\section{Tensor product of extremal weight crystals}
\subsection{A monoidal category of $\gl_\infty$-crystals}\label{Cat C}
Let $\mathcal{C}$ be the category of $\gl_{\infty}$-crystals, where each
object $B$ in $\mathcal{C}$ satisfies the following conditions;
\begin{itemize}
\item[($C1$)] there exists a finite subset $S\subset \cP\times \cP$ such that each connected component of $B$ is isomorphic to $\cB_{\mu,\nu}$ or $\cB_{\mu,\nu}\otimes \B(\Lambda_\lambda)$ for some $(\mu,\nu)\in S$ and $\lambda\in\Z_+^n$,

\item[($C2$)]  the number of connected components of $B$ isomorphic to $B(\Lambda)$ is finite for each $\Lambda\in P$,
\end{itemize}
and a morphism is a crystal morphism.

\begin{thm}\label{category C}
$\mathcal{C}$ is a monoidal category under tensor product of
crystals.
\end{thm}
\pf  It is enough to show that $B\otimes B'\in \mathcal{C}$ for
$B,B'$ in $\mathcal{C}$ since the map sending $(b_1\otimes
b_2)\otimes b_3$ to $b_1\otimes (b_2\otimes b_3)$ is an isomorphism
of crystals for $B_i\in \mathcal{C}$ and $b_i\in B_i$ ($i=1,2,3$).
By $(C1)$, it suffices to prove the case when
\begin{equation*}
\begin{split}
B&=\bigsqcup_{m\geq 1}\bigsqcup_{\lambda\in\Z_+^m}\cB_{\mu,\nu}\otimes\B(\Lambda_\lambda)^{\oplus c_{\lambda\mu\nu}}=
\cB_{\mu,\nu}\otimes\left(\bigsqcup_{m\geq 1}\bigsqcup_{\lambda\in\Z_+^m}\B(\Lambda_\lambda)^{\oplus c_{\lambda\mu\nu}}\right),\\
B'&=\bigsqcup_{n\geq 1}\bigsqcup_{\eta\in\Z_+^n}\cB_{\sigma,\tau}\otimes\B(\Lambda_\eta)^{\oplus c_{\eta \sigma\tau}}=
\cB_{\sigma,\tau}\otimes\left(\bigsqcup_{n\geq 1}\bigsqcup_{\eta\in\Z_+^n}\B(\Lambda_\eta)^{\oplus c_{\eta \sigma\tau}}\right)
\end{split}
\end{equation*}
for $\mu,\nu,\sigma,\tau\in\cP$  and
$c_{\lambda \mu \nu},c_{\eta \sigma  \tau}\in\Z_{\geq 0}$.

\textsc{Step 1}. Suppose that $\mu,\nu,\sigma,\tau=\emptyset$. Then
for $\zeta\in\Z_+^{l}$, the multiplicity of $\B(\Lambda_\zeta)$ in
$B\otimes B'$ is equal to $\sum_{\lambda,\eta}c_{\lambda\,\emptyset\,\emptyset }c_{\eta\,\emptyset\,\emptyset }
c^\zeta_{\lambda \eta}$ by Proposition \ref{LR HW}.  Since there are only finitely many $\lambda\in\Z_+^m$ and $\eta\in\Z_+^n$ such that $m+n=l$ and  $c^\zeta_{\lambda \eta}\neq 0$, it is a
well-defined integer.

\textsc{Step 2}. Let
$B''=\left(\bigsqcup_{m\geq 1}\bigsqcup_{\lambda\in\Z_+^m}\B(\Lambda_\lambda)^{\oplus
c_{\lambda\mu\nu}}\right)\otimes \cB_{\sigma,\tau}$.  By
Proposition \ref{commmutativity of BmBnvee} and Corollary
\ref{Pieri-2}, we have
\begin{equation*}
B''\simeq  
\bigsqcup_{(\alpha,\beta)\in S}\cB_{\alpha,\beta}\otimes\left(\bigsqcup_{m\geq 1}\bigsqcup_{\gamma\in\Z_+^{m}} \B(\Lambda_\gamma)^{\oplus d_{\gamma\alpha\beta}}\right)
\end{equation*}
for some finite subset $S$ of $\cP\times \cP$ and $d_{\gamma\alpha\beta}\in\Z_{\geq 0}$.
This implies that $B''\in\mathcal{C}$. By \textsc{Step 1}, we have $B'''=B''\otimes \left(\bigsqcup_{n\geq 1}\bigsqcup_{\eta\in\Z_+^n}\B(\Lambda_\eta)^{\oplus c_{\eta \sigma\tau}}\right)\in\mathcal{C}$. Finally, by Proposition \ref{LR for BmBn}, $B\otimes B'=\cB_{\mu,\nu}\otimes B'''\in\mathcal{C}$.
\qed\vskip 3mm

\begin{rem}{\rm
(1) A connected component of the tensor product $\B(\Lambda)\otimes
\B(\Lambda')$ with $\Lambda\in P_{-m}$ and $\Lambda'\in P_n$ for $m,n>0$ is not necessarily isomorphic to an extremal weight
crystal. For example, consider $B=\B(-\Lambda_0)\otimes
\B(\Lambda_0)$. We can check that any element in $B$
is connected to $u_{-\Lambda_0}\otimes u_{\Lambda_0}$, and hence
$B$ is connected. Suppose that $B$ is isomorphic to an extremal
weight crystal. By Corollary \ref{level 0 extremal}, $B\simeq
\B_{\mu,\nu}$ for some $\mu,\nu\in\cP$. By Remark \ref{remark on
bmn}, for given $S\otimes T\in \B_{\mu,\nu}$ and $p,q\in \Z$ with
$p<q$, there exist $i_1,\ldots,i_r\in \Z$ such that
$\td{x}_k\left(\te_{i_1}\cdots\te_{i_r} (S\otimes T) \right)={\bf
0}$ for $x=e,f$ and $k\in [p,q]$. Suppose that $S\otimes T$ is
equivalent to $u_{-\Lambda_0}\otimes u_{\Lambda_0}$. Then we have
\begin{equation*}
\te_{i_1}\cdots\te_{i_r} (u_{-\Lambda_0}\otimes u_{\Lambda_0})=
\left(\te_{i_1}\cdots\te_{i_r}  u_{-\Lambda_0}\right)\otimes u_{\Lambda_0}.
\end{equation*}
On the other hand, we have
\begin{equation*}
\begin{split}
&\tf_0\left(\te_{i_1}\cdots\te_{i_r} (u_{-\Lambda_0}\otimes u_{\Lambda_0})\right)\\
&=
\begin{cases}
\left(\te_{i_1}\cdots\te_{i_r}  u_{-\Lambda_0}\right)\otimes \tf_0 u_{\Lambda_0}, & \text{if $\varphi_0(\te_{i_1}\cdots\te_{i_r}  u_{-\Lambda_0})=0$}, \\
\left(\tf_0\te_{i_1}\cdots\te_{i_r}  u_{-\Lambda_0}\right)\otimes u_{\Lambda_0}, & \text{otherwise}, \\
\end{cases}
\end{split}
\end{equation*}
which is not ${\bf 0}$ in any case. This is a contradiction. Hence
$B$ is not isomorphic to an extremal weight crystal. In general, for
$\Lambda,\Lambda'\in P^+$, we can check by similar arguments that
$\B(-\Lambda)\otimes \B(\Lambda')$ is a connected regular crystal
but not isomorphic to an extremal weight crystal.

(2) When we consider a tensor product of arbitrary extremal weight
crystals, the multiplicity of each connected component is not
necessarily finite. For example, the multiplicity of $\B(\Lambda_0)$
in $\B(\Lambda_0)\otimes\B(\Lambda_0)\otimes\B(-\Lambda_0)$ is
infinite.}
\end{rem}

\subsection{Grothendieck ring}
Let $\K$ be the additive abelian group generated by the symbol
$[B]$ ($B\in\mathcal{C}$) subject to the relations; $[B]=[B']$ if
$B\simeq B'$ and $[B\sqcup B']=[B]+[B']$ for
$B, B'\in\mathcal{C}$.

\begin{prop}
$\K$ is an associative $\Z$-algebra with $1=[T_0]$ under the
multiplication $[B]\cdot[B']=[B\otimes B']$ for
$B,B'\in\mathcal{C}$.
\end{prop}

Let $\mathcal{C}_n$, $\mathcal{C}^{\rm h.w.}$ and
$\mathcal{C}^{\rm h.w.}_n$  ($n\geq 0$) be the full subcategories of $\mathcal{C}$
consisting of objects whose connected components are isomorphic to
$\B(\Lambda)$ for $\Lambda$ in $P_n$, $P^+$ and $P_n^+$,
respectively. We denote the corresponding subgroups of $\K$ by
$\K_{n}$, $\K^{\rm h.w.}$ and $\K^{\rm h.w.}_n$, respectively. Note that $\K_{0}$
is a subalgebra of $\K$.   By Proposition \ref{equi of bmnBLambda}, $[B]=[B']$
if and only if the multiplicities of each
$\B_{\mu,\nu}\otimes\B(\Lambda_\lambda)$ in $B$ and $B'$ are equal for $B,B'\in \mathcal{C}$. Hence, as $\Z$-modules, we have
\begin{equation}
\begin{split}
\K_n &\simeq \K_0\otimes_{\mathbb{Z}} \K^{\rm h.w.}_n, \\
\K&\simeq \K_0\otimes_{\mathbb{Z}} \K^{\rm h.w.},
\end{split}
\end{equation}
and $\K_n$, $\K$ are free $\K_0$-modules.

Let $\mathcal{C}^\vee$  be the category of $\gl_\infty$-crystals
with objects $B^\vee$ for $B\in\mathcal{C}$,  and let
$\mathcal{C}_{-n}$, $\mathcal{C}^{\rm l.w.}$ and
$\mathcal{C}^{\rm l.w.}_{-n}$  ($n\geq 0$) be its subcategories whose objects are $B^\vee$ for $B$ in
$\mathcal{C}_{n}$, $\mathcal{C}^{\rm h.w.}$ and $\mathcal{C}^{\rm
h.w.}_{n}$, respectively. We denote by $\K^\vee$, $\K_{-n}$,
$\K^{\rm l.w.}$ and $\K^{\rm l.w.}_{-n}$ the corresponding groups,
respectively. Then $\K^\vee$ is a $\Z$-algebra under tensor product
of crystals and isomorphic to $\K^{\rm opp}$, the opposite
$\Z$-algebra of $\K$. We have $\K_{-n}\simeq  \K^{\rm
l.w.}_{-n}\otimes_{\mathbb{Z}}\K_0$ and $\K^\vee\simeq \K^{\rm l.w.}\otimes_{\mathbb{Z}}\K_0$
as $\Z$-modules.

We denote by $Sym_x$ and $Sym_y$ the ring of symmetric functions in
$x$ and $y$, respectively, where $x=\{\,x_1,x_2,\ldots\,\}$ and
$y=\{\,y_1,y_2,\ldots\,\}$ are two sets of formal commuting
variables. Let $Sym_x\otimes_{\mathbb{Z}} Sym_y$ be the tensor product of
$\Z$-algebras with $\Z$-basis
$\{\,s_\mu(x)s_{\nu}(y)\,|\,\mu,\nu\in\cP\,\}$.

\begin{prop}\label{K0}
The assignment $s_\mu(x)s_\nu(y)\mapsto [\cB_{\mu,\nu}]$ $(\mu,\nu\in\cP)$ defines
an isomorphism  of $\Z$-algebras $\Phi : Sym_x\otimes_{\mathbb{Z}}
Sym_y\longrightarrow  \K_0.$
\end{prop}
\pf  Note that $\{\,[\cB_{\mu,\nu}]\,|\,\mu,\nu\in\cP\,\}$ is a
$\Z$-basis of $\K_0$  by  Proposition \ref{equi of bmn} and
Corollary \ref{level 0 extremal}. Let $\Phi : Sym_x\otimes_{\mathbb{Z}}
Sym_y\longrightarrow  \K_0$ be a  linear isomorphism sending $s_\mu(x)s_\nu(y)$ to $[\cB_{\mu,\nu}]$ $(\mu,\nu\in\cP)$.
It follows immediately from
Propositions \ref{commmutativity of BmBnvee} and \ref{LR for BmBn} that $\Phi$ is a homomorphism of algebras, and hence an isomorphism.
\qed\vskip 3mm

Let $\mathbb{Q}[P]$ be the group algebra of $P$ over $\mathbb{Q}$
with basis $\{\,e^\Lambda\,|\,\Lambda\in P\,\}$.  For $\Lambda\in
P^+$, let ${\rm ch}\B(\Lambda)=\sum_{b\in \B(\Lambda)}e^{{\rm
wt}(b)}$ be the character of $\B(\Lambda)$, which is a formal power
series in $\{\,e^{\Lambda_0}, e^{\pm\epsilon_i}\,  (i\in \Z)\,\}$ and
equal to the character of the irreducible highest weight
$U_q(\gl_\infty)$-module with highest weight $\Lambda$. Then $\{\,{\rm
ch}\B(\Lambda)\,|\,\Lambda\in P^+\,\}$ is linearly independent, and
${\rm ch}B$ is also well-defined for $B\in\mathcal{C}^{\rm h.w.}$
(use the formula for $e^{-n\Lambda_0}{\rm ch}\B(\Lambda_\lambda)$ ($\lambda\in\Z_+^n$) in \cite[Theorem 5.4]{K08} and then apply \cite[Proposition 3.18]{K08} for the well-definedness of ${\rm ch}B$ for $B\in \mathcal{C}^{\rm h.w.}$). Let ${R}$ be the $\Z$-algebra
spanned by $\{\,{\rm ch}B\,|\,B\in\mathcal{C}^{\rm h.w.}\,\}$. Then
the map $\psi : \K^{\rm h.w.} \rightarrow {R}$ given by
$\psi([B])={\rm ch}B$ is an algebra isomorphism.

Let $z=\{\,z_k\,|\,k\in\Z\,\}$ be another set of formal commuting
variables,  and let $\mathscr{R}$ be the ring of formal power series
in $z$ with coefficients in $\Z$. Then we have an
$A_\infty$-analogue of the fundamental theorem on symmetric
functions as follows.
\begin{prop}\label{K+}
The assignment $z_k\mapsto [\B(\pm\Lambda_k)]$ $(k\in\Z)$ defines isomorphisms of $\Z$-algebras
\begin{equation*}
\begin{split}
&\Psi_+ : \R \longrightarrow \K^{\rm h.w.}, \\ & \Psi_- : \R \longrightarrow \K^{\rm l.w.}, \\
\end{split}
\end{equation*}
respectively.
\end{prop}
\pf Let us identify $\K^{\rm h.w.}$ with ${R}$. Put $H_k={\rm
ch}\B(\Lambda_k)$ for $k\in \Z$.

For
$\lambda,\mu\in\mathbb{Z}_+^n$, we define $\lambda
> \mu$ if and only if there exists $i\geq 1$ such that
$\lambda_k=\mu_k$ for $1\leq k< i$ and $\lambda_i>\mu_i$. Then $>$
is a linear ordering on $\mathbb{Z}_+^n$.
Put
$H_{\mu}=\prod_{i=1}^nH_{\mu_i}$. Then we have
\begin{equation}\label{Kostka}
H_{\mu}=\sum_{\lambda\in\mathbb{Z}_+^n}K_{\lambda\,\mu}{\rm ch}\B(\Lambda_{\lambda})
\end{equation}
for some $K_{\lambda\,\mu}\in \Z_{\geq 0}$, where
$K_{\lambda\,\mu}=0$ unless $\lambda\geq \mu$ and
$\sum_i\lambda_i =\sum_i\mu_i$ (see (4.1) and (4.2) in \cite{K08}).
Note that we have a Jacobi-Trudi formula
\begin{equation}\label{jacobi}
{\rm ch}\B(\Lambda_\lambda)={\rm det}\left( H_{\lambda_i-i+j}\right)_{1\leq i,j\leq n}
\end{equation}
\cite{FKRW} (see also \cite[Theorem 4.5]{K08}).  

Define an algebra homomorphism
$\Psi_+ : \mathscr{R} \rightarrow R$  by $\Psi_+(z_k)=H_k$ for $k\in\Z$. Suppose
that $f(z)=c+\sum_{n\geq 1}\sum_{\mu\in\Z_+^n}c_\mu z_\mu$ is given (not necessarily
a finite sum), where $z_\mu=\prod_{i=1}^nz_{\mu_i}$. For $\lambda\in\Z_+^n$,  there are only finitely many $\mu\in\Z_+^n$ such that the multiplicity of $\B(\Lambda_\lambda)$
in $\B(\Lambda_{\mu_1})\otimes\cdots\otimes \B(\Lambda_{\mu_n})$ is non-zero by (\ref{Kostka}). Hence
$\Psi_+(f(z))=c+\sum_{n\geq 1}\sum_{\mu\in\Z_+^n}c_\mu H_\mu \in R$ is well-defined.

Given $c+\sum_{n\geq 1}\sum_{\lambda\in\Z_+^n}c'_\lambda  {\rm ch}\B(\Lambda_\lambda)\in R$, we have by (\ref{jacobi}) 
\begin{equation*}
c+\sum_{n\geq 1}\sum_{\lambda\in\Z_+^n}c'_\lambda  {\rm ch}\B(\Lambda_\lambda)=c+\sum_{n\geq 1}\sum_{\mu\in\Z_+^n}c_\mu H_\mu
\end{equation*}
for some $c_\mu\in\Z_{\geq 0}$. So $\Psi_+$ is surjective.
Also   by
(\ref{Kostka}), $\{\,H_\mu\,|\,\mu\in\Z_+^n,\ n\geq 1\,\}$ is linearly independent, which implies that  $\Psi_+$ is injective. Hence  $\Psi_+$ is an isomorphism. The proof for
$\Psi_-$ is almost the same. \qed\vskip 3mm

Let $t^\pm=\{\,t^\pm_1,t^\pm_2,t^\pm_3,\ldots\,\}$ be two sets of
mutually commuting formal variables.  Consider $$\cA=\R\,[t^+,t^-]$$
the free $\R$-module of the polynomials in $t^\pm$ over $\R$. For
$n\geq 1$, let $\cA_n=\R[t^\pm_1,\ldots,t^\pm_{n}]$. Note that as an $\R$-module
$$\cA_{n-1}\subset \cA_n=\cA_{n-1}[t^+_n,t^-_n], \ \ \
\cA=\sum_{n\geq 0}\cA_n,$$ where $\cA_0=\R$. Now, let us define a
$\Z$-algebra structure on $\cA$ inductively as follows;
\begin{itemize}
\item[(1)] $\cA_0=\R$ with the usual multiplication.

\item[(2)] Suppose that a $\Z$-algebra structure on $\cA_{n-1}$ is defined.
Then we define a multiplication on $\cA_n$  by
\begin{equation}\label{commutator}
t^\pm_{n} a= a t^\pm_n + \delta^\pm_n(a)\ \ \ (a\in\cA_{n-1}),
\end{equation}
    where $\delta^\pm_n$ is the derivations on $\cA_{n-1}$
given by
\begin{equation}\label{derivation delta}
\begin{split}
&\delta^\pm_n(t^\pm_{k})=0\ \ (1\leq k\leq n-1), \\
&\delta^\pm_n(z_k)=z_{k\mp 1}t^\pm_{n- 1}+ z_{k\mp 2} t^\pm_{n- 2}+
\cdots +z_{k\mp n}\ \ (k\in\Z)
\end{split}
\end{equation}
(that is, $\cA_{n}$ is an Ore extension \cite{Ore} of $\cA_{n-1}$
associated with derivations $\delta^\pm_n$).
\end{itemize}
Note that $\delta^\pm_n$ is well-defined on $\cA_{n-1}$ since for a generalized
partition $\mu$,   $\delta^\pm_n(z_{\mu})$ is a polynomial in
$t^\pm_1,\ldots,t^\pm_{n-1}$ having polynomials in $z$ as coefficients.

\begin{prop}\label{K}
The assignment $z_k\mapsto [\B(-\Lambda_k)]$,
$t^+_n\mapsto [\cB_{(1^n)}]$ and $t^-_n\mapsto[\cB_{(1^n)}^\vee]$  $(k\in \Z, n\geq 1)$
defines an isomorphism of $\Z$-algebras  $\Psi : \cA\longrightarrow
\K^\vee.$
\end{prop}
\pf Since the map sending $t^+_k t^-_l$ to $e_k(x)e_l(y)$ ($k,l\geq 1$) gives a
$\Z$-algebra isomorphism from $\Z[t^+,t^-]$ to $Sym_x\otimes_{\mathbb{Z}} Sym_y$,
composing with $\Phi$ in Proposition \ref{K0}  we have a
$\Z$-algebra isomorphism
\begin{equation*}\label{Psi0}
\Psi_0 : \Z[t^+,t^-]\longrightarrow \K_0
\end{equation*}
given  by $\Psi_0(t^+_n)=[\cB_{(1^n)}]$ and
$\Psi_0(t^-_n)=[\cB_{(1^n)}^\vee]$ for $n\geq 1$.

Since $\cA\simeq \R\otimes_{\mathbb{Z}} \Z[t^+,t^-]$ and $\K^\vee\simeq
\K^{\rm l.w.}\otimes_{\mathbb{Z}}\K_0$ as $\Z$-modules,  we have a $\Z$-linear
isomorphism
\begin{equation*}
\Psi=\Psi_-\otimes\Psi_0 : \cA  \longrightarrow  \K^\vee.
\end{equation*}

Next, we claim that $\Psi$ is an algebra homomorphism. For $n\geq
0$, let $\R_n=\Psi(\cA_n)$. We use induction on $n$ to  show that
$\Psi|_{\cA_n} : \cA_n \longrightarrow \R_n$ is an algebra
homomorphism. If $n=0$, then $\Psi|_{\cA_0}=\Psi_- : \R
\longrightarrow \K^{\rm l.w.}$ is an algebra isomorphism by
Proposition \ref{K+}. Suppose that $\Psi|_{\cA_{n}}$ is an algebra
homomorphism. Note that $\cA_{n+1}=\cA_n[t^+_{n+1},t^-_{n+1}]$. By
Lemmas \ref{babb commutes} and \ref{commuting relation-1}, we have
\begin{equation*}
\Psi|_{\cA_{n+1}}(t^\pm_{n+1} a - a t^\pm_{n+1} -
\delta^\pm_{n+1}(a))=0,
\end{equation*}
when $a=z_k$ ($k\in\Z$) and $a=t^\pm_k$ ($1\leq k\leq n$). 
It follows that $\Psi|_{\cA_{n+1}}$ preserves the multiplication,
and hence it is an algebra homomorphism. This completes the
induction. Therefore, $\Psi$ is an algebra isomorphism.
\qed\vskip 3mm

Let $s^\pm=\{\,s^\pm_1,s^\pm_2,s^\pm_3,\ldots,\}$ be two sets of
mutually commuting formal variables. Consider
$$\cD=\R_{\mathbb{Q}}\,[s^+,s^-]$$ the free $\R_{\mathbb{Q}}$-module of the polynomials in
$s_\pm$ over $\R_{\mathbb{Q}}={\mathbb{Q}}\otimes_{\mathbb{Z}} \R$. Then we define a $\mathbb{Q}$-algebra structure on $\cD$
by
\begin{equation}\label{commutator-2}
\begin{split}
&s^\pm_{n} z_k= z_k s^\pm_n + (-1)^{n-1}z_{k\mp n}, \ \ (n\geq 1, k\in\Z). \\
\end{split}
\end{equation}
One may regard $\cD$ as an Ore extension of $\R_{\mathbb{Q}}$ associated with
derivations $\gamma^\pm_n=(-1)^{n-1}\sum_{k\in\Z}z_{k\mp
n}\frac{\partial}{\partial z_k}$ ($n\geq 1$).

Let $u$ be a formal  variable. Put
\begin{equation}
F(u)= \sum_{k\in\Z}z_k u^k, \ \ E^\pm(u)=\exp\left(\sum_{n\geq
1}\frac{(-1)^{n-1}}{n}s^\pm_n u^n \right) \in \cD[[u]].
\end{equation}
Now, we
obtain the following characterization of
$\K^\vee_{\mathbb{Q}}=\mathbb{Q}\otimes_{\mathbb{Z}}\K^\vee$ (hence $\K_{\mathbb{Q}}=\mathbb{Q}\otimes_{\mathbb{Z}}\K$), which is
the main result in this section.
\begin{thm}\label{K-main}
There exists a $\mathbb{Q}$-algebra isomorphism $\Theta : \cD \longrightarrow \K^\vee_{\mathbb{Q}}$ such that
\begin{equation*}
\begin{split}
& \Theta\left(F(u)\right)=\sum_{k\in \Z}[\B(-\Lambda_k)]u^k,\\
& \Theta \left(E^+(u) \right)  = \sum_{n\geq 0}[\B_{(1^n)}]u^n,  \ \
  \Theta \left(E^-(u) \right)  = \sum_{n\geq 0}[\B_{(1^n)}^\vee]u^n .
\end{split}
\end{equation*}
\end{thm}
\pf For $n\geq 1$, set
\begin{equation*}
\hat{s}^+_n=\Psi_0^{-1}\circ\Phi(p_n(x)), \ \hat{s}^{\,-}_n=\Psi_0^{-1}\circ\Phi(p_n(y)).
\end{equation*}
Note that $\hat{s}^\pm=\{\,\hat{s}^\pm_n\,|\,n\geq 1\,\}$ is algebraically
independent over $\mathbb{Q}$, and $\cA_{\mathbb{Q}}=\mathbb{Q}\otimes_{\mathbb{Z}}\cA=\R_{\mathbb{Q}}\,[\hat{s}^+,\hat{s}^-]$.

First, we use induction on $n$ to show that for $k\in\Z$
\begin{equation}\label{snzk commuting rel}
\hat{s}^{\,\pm}_n z_k =z_k \hat{s}^{\pm}_n +(-1)^{n-1}z_{k\mp
n}.
\end{equation}
It is clear when $n=1$. Recall that
\begin{equation}\label{ep rel}
(-1)^{n-1}\hat{s}^{\pm}_n=nt^\pm_n
-\sum_{r=1}^{n-1}(-1)^{r-1}\hat{s}^{\pm}_r t^\pm_{n-r}.
\end{equation}
By induction hypothesis, we have for $k\geq 1$,
\begin{equation*}
\begin{split}
(-1)^{n-1}\hat{s}^{\pm}_n z_k=& \ nt^\pm_n
z_k-\sum_{r=1}^{n-1}(-1)^{r-1}\hat{s}^{\pm}_r t^\pm_{n-r}
z_k \\
=& \ n\left( z_kt^\pm_n +z_{k\mp 1}t^\pm_{n-1}+\cdots+ z_{k\mp
n}\right)\\ &-\sum_{r=1}^{n-1}(-1)^{r-1}\sum_{j=0}^{n-r}(z_{k\mp
j}\hat{s}^{\pm}_r+(-1)^{r-1}z_{k\mp r\mp j})t^{\pm}_{n-r-j}. \\
\end{split}
\end{equation*}
Using (\ref{ep rel}), it is straightforward to check that the
(right) coefficient of $z_{k\mp i}$ in the last equation is
\begin{equation*}
\begin{cases}
(-1)^{n-1} \hat{s}^\pm_n, & \text{if $i=0$}, \\
1, & \text{if $i=n$}, \\
0, & \text{otherwise}.
\end{cases}
\end{equation*}
This proves (\ref{snzk commuting rel}) and completes the induction.

By (\ref{commutator-2}) and (\ref{snzk commuting rel}), we obtain a
$\mathbb{Q}$-algebra isomorphism $\theta : \cD \longrightarrow \cA_{\mathbb{Q}}$ such
that $\theta(t_k)=t_k$ and $\theta(s^\pm_n)=\hat{s}^\pm_n$ for $k\in\Z$ and $n\geq 1$. Since
$$\exp\left(\sum_{n\geq 1}\frac{(-1)^{n-1}}{n}\hat{s}^\pm_n u^n \right)=\sum_{n\geq 0}t^\pm_n u^n$$
(cf.\cite{Mac95}), we obtain the required isomorphism $\Theta=\Psi\circ \theta$
by Proposition \ref{K}. \qed

\begin{cor}
$\K_{\mathbb{Q}}$ is isomorphic to $\cD^{\rm opp}$ as a $\mathbb{Q}$-algebra, where $\mathscr{D}^{\rm \, opp}$ denotes the opposite algebra of
$\mathscr{D}$.
\end{cor}

Note that there exists an involution $\omega : \cD \longrightarrow \cD$ determined by
\begin{equation}\label{omega}
\omega(z_k)=z_{-k},\ \ \ \omega(s^\pm_{n})=s^\mp_n
\end{equation}
for $k\in\Z$ and $n\geq 1$, which is well-defined by
(\ref{commutator-2}).

\begin{prop}
For $\lambda\in \Z_+^n$ and $\mu,\nu\in\cP$, we have
$$\omega\left( [\B(-\Lambda_\lambda)\otimes\cB_{\mu,\nu}]\right)=[\B(-\Lambda_{\lambda^*})\otimes\cB_{\nu,\mu}],$$
where $\lambda^*=(-\lambda_n,\ldots,-\lambda_1)$.
\end{prop}
\pf It follows directly from (\ref{jacobi}) that $\omega
[\B(-\Lambda_\lambda)]=[\B(-\Lambda_{\lambda^*})]$. Since
$\Psi^{-1}_0\circ\Phi(p_m(x))=\hat{s}^+_m$ and
$\Psi^{-1}_0\circ\Phi(p_m(y))=\hat{s}^-_m$ for $m\geq 1$, we have $$\left(\Phi^{-1}\circ \omega \circ
\Phi\right)(s_{\mu}(x)s_{\nu}(y))=s_{\nu}(x)s_{\mu}(y),$$ which
implies $\omega([\B_{\mu,\nu}])=[\B_{\nu,\mu}]$. Since $\omega$ is
an algebra homomorphism, we obtain the above identity. \qed

\subsection{Littlewood-Richardson rule}

Let $u$ be a formal variable. Put
\begin{equation}
\begin{split}
&\mathscr{F}(u)= \sum_{k\in \Z}[\B(\Lambda_k)]u^k,\ \ \mathscr{E}(u)=\sum_{n\geq 0}[\B_{(1^n)}]u^n, \ \  \mathscr{E}^\vee(u)=\sum_{n\geq 0}[\B^\vee_{(1^n)}]u^n.
\end{split}
\end{equation}

\begin{lem}\label{Cauchy identities}
For $n\geq 1$, we have {\rm
\begin{equation*}
\begin{split}
\mathscr{F}(x_{[n]})&=\sum_{\lambda\in\Z_+^n} [\B(\Lambda_\lambda)]s_\lambda(x_{[n]}), \\
\mathscr{E}(x_{[n]})&=\sum_{\substack{\mu \in \cP \\ \mu_1\leq n}} [\cB_{\mu}] s_{\mu'}(x_{[n]}), \ \
\mathscr{E}^\vee(x_{[n]}) =\sum_{\substack{\nu \in \cP \\ \nu_1\leq n}} [\cB_{\nu}^\vee] s_{\nu'}(x_{[n]}),
\end{split}
\end{equation*}}
where {\rm $\mathscr{F}(x_{[n]})=\prod_{k=1}^n\mathscr{F}(x_k)$,
$\mathscr{E}(x_{[n]})=\prod_{k=1}^n\mathscr{E}(x_k)$} and {\rm $\mathscr{E}^\vee(x_{[n]})=\prod_{k=1}^n\mathscr{E}^\vee(x_k). $}
\end{lem}
\pf Consider the $(\gl_\infty,\gl_n)$-character associated with the
decomposition  in Proposition \ref{duality-1}. Then we obtain the first
identity   by replacing ${\rm ch}\B(\Lambda_\lambda)$ with
$[\B(\Lambda_\lambda)]$. The other two identities are obtained by
considering the $(\gl_\infty,\gl_n)$-character associated with the
decomposition in Proposition \ref{duality on En} and using the
isomorphism in Proposition \ref{K0}. \qed\vskip 3mm

Let $u,v$ be commuting formal variables. By Lemma \ref{commuting relation-1}, we have
\begin{equation}\label{commuting relation-2}
\begin{split}
\mathscr{F}(u) \mathscr{E}(v)&= \mathscr{E}(v)\mathscr{F}(u)\frac{1}{(1-u^{-1}v)}, \\
\mathscr{F}(u) \mathscr{E}^\vee(v)&= \mathscr{E}^\vee(v)\mathscr{F}(u)\frac{1}{(1-uv)}.
\end{split}
\end{equation}
Applying (\ref{commuting relation-2}) successively, we obtain the following identities.
\begin{lem}\label{Cauchy-noncommutative} For $m,n\geq 1$,
we have {\rm
\begin{equation*}
\begin{split}
\mathscr{F}(x_{[m]})\mathscr{E}(y_{[n]})&=\mathscr{E}(y_{[n]})\mathscr{F}(x_{[m]})\frac{1}{\prod_{i\in[m],j\in[n]}(1-x_i^{-1}y_j)}, \\
\mathscr{F}(x_{[m]})\mathscr{E}^\vee(y_{[n]})&=\mathscr{E}^\vee(y_{[n]})\mathscr{F}(x_{[m]})\frac{1}{\prod_{i\in[m],j\in[n]}(1-x_iy_j)}.
\end{split}
\end{equation*}}
\end{lem}
\vskip 3mm

\begin{lem}\label{Pieri-general} Let $\lambda,\eta\in\Z_+^m$ and $\mu,\nu\in\cP$ be given.
\begin{itemize}
\item[(1)] The multiplicity of $\cB_\nu\otimes \B(\Lambda_\eta)$ in
$\B(\Lambda_\lambda)\otimes \cB_\mu$ is
$$\sum_{\substack{\gamma\in\cP \\ \ell(\gamma)\leq m}}
c^{\lambda}_{\eta\, \gamma^\ast}c^{\mu'}_{\nu' \gamma},$$
where $\gamma^\ast=(\ldots,-\gamma_2,-\gamma_1)\in\Z_+^m$.

\item[(2)] The multiplicity of $\cB^\vee_\nu\otimes \B(\Lambda_\eta)$ in  $\B(\Lambda_\lambda)\otimes \cB^\vee_\mu$ is
    $$\sum_{\substack{\gamma\in\cP \\ \ell(\gamma)\leq m}} c^{\lambda}_{\eta \gamma}c^{\mu'}_{\nu' \gamma}.$$
\end{itemize}
\end{lem}
\pf (1) Choose $n\geq 1$ such that $\ell(\mu')\leq n$. The left-hand side of the first identity in Lemma \ref{Cauchy-noncommutative} is given by
\begin{equation}\label{LHS}
\sum_{\lambda\in\Z_+^m}\sum_{\substack{\mu\in\cP \\
\mu_1\leq
n}}[\B(\Lambda_\lambda)\otimes\cB_\mu]s_{\lambda}(x_{[m]})s_{\mu'}(y_{[n]}).
\end{equation}
On the other hand, the right-hand side is
\begin{equation}\label{RHS}
\begin{split}
&\sum_{\eta\in\Z_+^m}\sum_{\substack{\nu\in\cP \\
\nu_1\leq
n}}[\cB_\nu\otimes\B(\Lambda_\eta)]s_{\eta}(x_{[m]})s_{\nu'}(y_{[n]})
\sum_{\substack{\gamma\in\cP \\ \ell(\gamma)\leq m,n}}s_{\gamma}(x^{-1}_{[m]})s_{\gamma}(y_{[n]}) \\
&=\sum_{\lambda\in\Z^m_+}\sum_{\substack{\mu\in\cP \\ \mu'\leq
n}}\sum_{\eta,\nu,\gamma} c^{\lambda}_{\eta\, \gamma^\ast}c^{\mu'}_{\nu'
\gamma}[\cB_{\nu}\otimes\B(\Lambda_\eta)]
s_{\lambda}(x_{[m]})s_{\mu'}(y_{[n]}).
\end{split}
\end{equation}
Since
$\{\,s_{\lambda}(x_{[m]})s_{\mu'}(y_{[n]})\,|\,\lambda\in\Z_+^m,
\mu\in\cP \text{ with $\ell(\mu')\leq n$}\,\}$  is linearly
independent, we have
\begin{equation*}
[\B(\Lambda_\lambda)\otimes\cB_\mu]=\sum_{\eta,\nu}\left(\sum_{\gamma}
c^{\lambda}_{\eta\, \gamma^\ast}c^{\mu'}_{\nu'
\gamma}\right)[\cB_{\nu}\otimes\B(\Lambda_\eta)]\in \K
\end{equation*}
by comparing (\ref{LHS}) and (\ref{RHS}).  Hence we obtain the
required multiplicity.

(2) The proof is almost the same as in (1).  We leave the details to
the reader. \qed\vskip 3mm

Combining Lemma \ref{Pieri-general} (1) and (2),  we have the
following decomposition of the tensor product of a highest weight
crystal and a level zero extremal weight crystal.

\begin{prop}\label{commuting relation-3}
For $\lambda\in\Z_+^m$ and $\mu,\nu\in \cP$,
we have
\begin{equation*}
\B(\Lambda_{\lambda})\otimes \cB_{\mu,\nu}
\simeq\bigsqcup_{\substack{\rho\in\Z_+^m \\ \sigma,\tau\in\cP}}
\cB_{\sigma,\tau}\otimes\B(\Lambda_{\rho})^{\oplus
c^{(\lambda,\mu,\nu)}_{(\rho,\sigma,\tau)}},
\end{equation*}
where
\begin{equation*}
c^{(\lambda,\mu,\nu)}_{(\rho,\sigma,\tau)}=\sum_{\eta\in\Z_+^m}\sum_{\substack{\alpha,\beta\in
\cP \\ \ell(\alpha),\ell(\beta)\leq m}}
c^{\lambda}_{\eta\, \alpha^\ast}c^{\mu'}_{\sigma'\,\alpha}c^{\eta}_{\rho\,\beta}c^{\nu'}_{\tau'\,\beta}.
\end{equation*}
\end{prop}

Now, we can describe the Littlewood-Richardson rule of extremal
weight  crystals of non-negative level, which is the main result
in this paper.

\begin{thm}\label{Extremal LR rule}
For $\lambda\in\Z_+^m$, $\rho\in\Z_+^n$,
and  $\mu,\nu,\sigma,\tau\in\cP$, we have
\begin{equation*}
\left(\cB_{\mu,\nu}\otimes
\B(\Lambda_{\lambda})\right)\otimes \left(\cB_{\sigma,\tau}\otimes
\B(\Lambda_{\rho})\right)\simeq \bigsqcup_{\substack{\zeta\in\Z_+^{m+n} \\ \eta,\theta\in\cP}} \cB_{\eta,\theta}\otimes
\B(\Lambda_\zeta)^{\oplus c^{(\zeta,\eta,\theta)}_{(\lambda,\mu,\nu),(\rho,\sigma,\tau)}},
\end{equation*}
where
\begin{equation*}
c^{(\zeta,\eta,\theta)}_{(\lambda,\mu,\nu),(\rho,\sigma,\tau)}=
\sum_{\alpha\in\Z_+^m }\sum_{\beta,\gamma\in \cP }c^{\zeta}_{\alpha
\rho}c^{\eta}_{\beta \mu}c^{\theta}_{\gamma
\nu}c^{(\lambda,\sigma,\tau)}_{(\alpha,\beta,\gamma)}
\end{equation*}
and $c^{(\lambda,\sigma,\tau)}_{(\alpha,\beta,\gamma)}$ is defined
in Proposition \ref{commuting relation-3}.
\end{thm}
\pf It follows from Propositions \ref{commmutativity of BmBnvee},
\ref{LR for BmBn}, \ref{LR HW}, and \ref{commuting
relation-3}. \qed

\begin{rem}{\rm
We have the same Littlewood-Richardson rule for the crystals in
$\mathcal{C}^\vee$ by taking the dual of the decomposition in
Theorem \ref{Extremal LR rule}. }
\end{rem}

\section{Differential operators on lowest weight character ring}

\subsection{A twisted action of $\K^\vee$ on $\K^{\rm l.w.}$ }

We define for $B\in \mathcal{C}^\vee$
\begin{equation*}
{\rm pr}(B)=\{\,b\in B\,|\, \text{$\exists\ r\geq 1$ such that
$\tf_{i_1}\cdots\tf_{i_r}b={\bf 0}$  for all $i_1,\ldots,i_r\in\Z$}
\,\}.
\end{equation*}
We can check that ${\rm pr}(B)$ is a union of connected components
$B'$ of $B$ such that $B' \in \mathcal{C}^{\rm l.w.}$. Hence ${\rm pr}$ is a
functor from $\mathcal{C}^\vee$ to $\mathcal{C}^{\rm l.w.}$, and by definition ${\rm pr}(B\sqcup B')\simeq {\rm pr}(B)\,\sqcup\, {\rm pr}(B')$
for $B, B'\in\mathcal{C}^\vee$. Consider the 
composite of following two functors
\begin{equation}
\begin{array}{ccccc}
  \mathcal{C}^\vee \times \mathcal{C}^{\rm l.w.} & \stackrel{\otimes}{\longrightarrow} & \mathcal{C}^\vee & \stackrel{{\rm
pr}}{\longrightarrow} &  \mathcal{C}^{\rm l.w.} \\
  (B,B') & \longmapsto & B\otimes B' & \longmapsto & {\rm pr}(B\otimes
  B')
\end{array}\ .
\end{equation}
Then it induces a $\Z$-algebra homomorphism (or $\K^\vee$-module structure on $\K^{\rm l.w.}$)
\begin{equation}\label{rho}
\rho : \K^\vee \longrightarrow {\rm End}_\Z(\K^{\rm l.w.}),
\end{equation}
where $\rho\left([B]\right)\left([B']\right) = [{\rm pr}(B\otimes B')]$
for $[B]\in \K^\vee$ and $[B']\in \K^{\rm l.w.}$. Hence $\K^{\rm
l.w.}$ is a left $\K^\vee$-module. Moreover, when restricted to
$\K_0$, each $\K^{\rm l.w.}_{-n}$ ($n\geq 0$) is a $\K_0$-submodule of $\K^{\rm
l.w.}$. Since the action of $\K^{\rm l.w.}$ on $\K^{\rm l.w.}$ is
nothing but the left multiplication in $\K^\vee$, we will focus on
the $\K_0$-module structure on $\K^{\rm l.w.}$.

Recall that by Proposition \ref{K}, we may identify $\K^\vee$ with
$\cA=\R\,[t^+,t^-]$, while $\K^{\rm l.w.}$ and $\K_0$ are identified
with $\R$ and $\Z[t^+,t^-]$, respectively. Let   $\R_n$ be the subspace of $\R$ consisting of formal
power series in $z$ of degree $n$. Then $\K^{\rm l.w.}_{-n}$
corresponds to $\R_n$.

With these identification, the left $\K^\vee$-module $\K^{\rm l.w.}$ corresponds to a left $\cA$-module
\begin{equation*}
\cA/\mathscr{I}, \ \ \ \ \mathscr{I}=\sum_{m\geq 1}\left(\cA t^+_m +
\cA t^-_m \right),
\end{equation*}
which can be identified with $\R$ as a  $\Z$-module.  We still denote
this $\cA$-module structure on $\R$ by $\rho : \cA \longrightarrow
{\rm End}_\Z(\R)$. Note that the action of $\R\subset \cA$ is the
usual multiplication on $\R$ and $\R_n$ is a $\Z[t^+,t^-]$-submodule
of $\R$.

\subsection{The action of $\K_0$ on $\K^{\rm l.w.}$}
Let us describe the action of $\Z[t^+,t^-]$ on $\R$ more explicitly.
Recall the following correspondences (see Propositions \ref{K+} and
\ref{K});
\begin{equation}\label{correspondence}
\begin{array}{ccccc}
  Sym_x\otimes Sym_y & \stackrel{\Phi}{\longrightarrow} & \K_0 & \stackrel{\Psi_0}{\longleftarrow} & \Z[t_+,t_-]  \\
  s_\mu(x)s_\nu(y) & \longleftrightarrow & [\cB_{\mu,\nu}] & \longleftrightarrow & t^+_{\{\mu\}}t^-_{\{\nu\}} \\
  e_m(x)e_n(y) & \longleftrightarrow & [\cB_{(1^m),(1^n)}] & \longleftrightarrow &
  t^+_mt^-_n
\end{array}
\end{equation}
Here $t^\pm_{\{\mu\}}={\rm det}(t^\pm_{\mu'_i-i+j})_{1\leq i,j\leq
\ell(\mu')}$, where we assume $t^\pm_0=1$ and $t^\pm_k=0$ for $k<0$.

For $n\geq 1$, let
\begin{equation}
\begin{split}
\textsf{p}^\pm_n&=\rho\left(\hat{s}^\pm_n\right),
\end{split}
\end{equation}
where $\hat{s}^+_n=\Psi^{-1}_0\circ\Phi(p_n(x))$ and $\hat{s}^-_n=\Psi^{-1}_0\circ\Phi(p_n(y))$. The following is an immediate
consequence of Theorem \ref{K-main} (see (\ref{snzk commuting rel})).
\begin{prop}\label{p action} For $n\geq 1$,
{\rm $\textsf{p}^\pm_n=\gamma^\pm_n=(-1)^{n-1}\sum_{k\in\Z}z_{k\mp
n}\frac{\partial}{\partial z_k}$.}
\end{prop}

For $\lambda,\mu\in\Z_+^n$, put
\begin{equation}
z_{\{\lambda/\mu\}} ={\rm det}(z_{\lambda_i-\mu_j-i+j})_{1\leq
i,j\leq n}.
\end{equation}
By Proposition \ref{K+}, each element $f(z)$ in $\R$ can be written
uniquely as $f(z)=c_0+\sum_{n\geq 1}\sum_{\lambda\in\Z_+^n}c_\lambda
z_{\{\lambda\}}$ for some $c_0, c_\lambda\in\Z$. 

\begin{lem}\label{Jacobi for Ainfty} For $\lambda,\mu\in\Z_+^n$, we have
$$z_{\{\lambda/\mu\}}=\sum_{\nu\in\Z_+^n}c^{\lambda}_{\mu \nu}z_{\{\nu\}}.$$
\end{lem}
\pf  For $p\ll \min\{0,\lambda_n\}$,
let
$$\B(\Lambda_\lambda)_{>p}=\{\,\tf_{i_1}\cdots\tf_{i_r}u_{\Lambda_\lambda}\,|\,r\geq
0,\ i_1,\ldots,i_r\in [p+1,\infty)\,\}\subset \B(\Lambda_\lambda).$$
Note that $\B(\Lambda_\lambda)_{>p}$ is the connected
$\gl_{[p+1,\infty)}$-subcrystal of $\B(\Lambda_\lambda)$ including
$u_{\Lambda_\lambda}$. Put $D_{p,n}=\prod_{p+1\leq i\leq
0}x_i^{-n}$. Then it is easy to see that
\begin{equation*}
e^{-n\Lambda_0}{\rm
ch}\B(\Lambda_\lambda)_{>p}=D_{p,n}s_{(\lambda-(p^n))'}(x_{[p+1,\infty)}),
\end{equation*}
where $x_{i}=e^{\epsilon_{i}}$ for $i\in [p+1,\infty)$. For $k\in [p+1,\infty)$, put
\begin{equation*}
\widehat{e}_k(x_{[p+1,\infty)})=D_{p,1}e_{k-p}(x_{[p+1,\infty)}).
\end{equation*}
Then $e^{-\Lambda_0}{\rm
ch}\B(\Lambda_k)_{>p}=\widehat{e}_k(x_{[p+1,\infty)})$ and 
$\widehat{e}_k(x_{[p+1,\infty)})$ has a well-defined limit when
$p\rightarrow -\infty$, which is equal to $e^{-\Lambda_0}{\rm
ch}\B(\Lambda_k)=e^{-\Lambda_0}H_k$ (cf. \cite[Section
3.3]{K08-JACO}).

We have
\begin{equation*}
\begin{split}
&{\rm det}(\widehat{e}_{\lambda_i-\mu_j-i+j}(x_{[p+1,\infty)}))_{1\leq i,j\leq n}\\&=D_{p,n}
{\rm det}({e}_{\lambda_i-\mu_j-i+j-p}(x_{[p+1,\infty)}))_{1\leq i,j\leq n} \\
&=D_{p,n}
s_{\left(\lambda-(p^n)-(\mu_n^n)\right)'/\left(\mu-(\mu_n^n)\right)'}(x_{[p+1,\infty)}).
\end{split}
\end{equation*}
On the other hand, we have
\begin{align*}
&s_{\left(\lambda-(p^n)-(\mu_n^n)\right)'/\left(\mu-(\mu_n^n)\right)'}(x_{[p+1,\infty)})) \\
&=\sum_{\nu\in\Z_+^n}c^{\left(\lambda-(p^n)-(\mu_n^n)\right)'}_{\left(\mu-(\mu_n^n)\right)' \left(\nu-(p^n) \right)'}s_{\left(\nu-(p^n) \right)'}(x_{[p+1,\infty)}) \\
&=\sum_{\nu\in\Z_+^n}c^{\lambda-(p^n)-(\mu_n^n)}_{\mu-(\mu_n^n)\, \nu-(p^n)  }s_{\left(\nu-(p^n) \right)'}(x_{[p+1,\infty)}) \\
&=\sum_{\nu\in\Z_+^n}c^{\lambda}_{\mu \nu}s_{\left(\nu-(p^n) \right)'}(x_{[p+1,\infty)}).
\end{align*}
Therefore, we have
\begin{equation}\label{skewJacobi}
\begin{split}
&{\rm det}(\widehat{e}_{\lambda_i-\mu_j-i+j}(x_{[p+1,\infty)}))_{1\leq i,j\leq n}\\
&=D_{p,n}\sum_{\nu\in\Z_+^n}c^{\lambda}_{\mu \nu}s_{\left(\nu-(p^n) \right)'}(x_{[p+1,\infty)}) \\
&=D_{p,n}\sum_{\nu\in\Z_+^n}c^{\lambda}_{\mu \nu}{\rm det}({e}_{\nu_i-i+j-p}(x_{[p+1,\infty)}))_{1\leq i,j\leq n} \\
&=\sum_{\nu\in\Z_+^n}c^{\lambda}_{\mu \nu}{\rm
det}(\widehat{e}_{\nu_i-i+j}(x_{[p+1,\infty)}))_{1\leq i,j\leq n}.
\end{split}
\end{equation}
Taking $p\rightarrow -\infty$ and then multiplying $e^{n\Lambda_0}$ on both sides of (\ref{skewJacobi}),
we have
\begin{equation*}
{\rm det}(H_{\lambda_i-\mu_j-i+j})_{1\leq i,j\leq
n}=\sum_{\nu\in\Z_+^n}c^{\lambda}_{\mu \nu}{\rm
det}(H_{\nu_i-i+j})_{1\leq i,j\leq n}.
\end{equation*}
By Proposition \ref{K+}, we have $z_{\{\lambda/\mu\}}=\sum_{\nu\in\Z_+^n}c^{\lambda}_{\mu \nu}z_{\{\nu\}}.$
\qed\vskip 3mm

\begin{rem}{\rm
We may realize $z_{\{\lambda/\mu \}}$ or
$\Psi_{\pm}(z_{\{\lambda/\mu \}})$ as a weight generating function
for certain pairs of semistandard tableaux given in \cite[Definition
4.9]{K08}. A bijective proof of Lemma \ref{Jacobi for Ainfty} using
this realization is also given in \cite[Theorem 4.11]{K08}. }
\end{rem}

For $\mu\in\cP$, we put
\begin{equation}
\begin{split}
& \textsf{s}^\pm_\mu= \rho\left(t^\pm_{\{\mu\}}\right).
\end{split}
\end{equation}

\begin{thm}\label{rho s lambda} For $\lambda\in\Z_+^n$ and $\mu\in\cP$, we have
$$\text{\rm $\textsf{s}^+_{\mu'}\left(z_{\{\lambda\}}\right)$}=
\begin{cases}
z_{\{\lambda/\mu\}}, & \text{if $\ell(\mu)\leq n$}, \\
0, & \text{otherwise},
\end{cases},\ \ \
\text{\rm $\textsf{s}^-_{\mu'}\left(z_{\{\lambda\}}\right)$}=
\begin{cases}
z_{\{\lambda/\mu^\ast\}}, & \text{if $\ell(\mu)\leq n$}, \\
0, & \text{otherwise},
\end{cases}$$
where $\mu^*=(\ldots,-\mu_2,-\mu_1)\in\Z_+^n$.
\end{thm}
\pf Note that the coefficient of $z_{\{\nu\}}$
for $\nu\in\Z_+^n$ in
$\textsf{s}^+_{\mu'}\left(z_{\{\lambda\}}\right)$ is equal to the
multiplicity of $\B(-\Lambda_\nu)$ in $\cB_{\mu'}\otimes
\B(-\Lambda_\lambda)$, or that of $\B(\Lambda_\nu)$ in
$\B(\Lambda_\lambda)\otimes \cB_{\mu'}^\vee$. By Lemma
\ref{Pieri-general}, it is equal to
\begin{equation}\label{multiplicity}
\sum_{\substack{\gamma\in\cP \\ \ell(\gamma)\leq n}}
c^{\lambda}_{\nu\, \gamma}c^{\mu}_{\emptyset \gamma}.
\end{equation}
Since
$$c^{\mu}_{\emptyset \gamma}=
\begin{cases}
1, & \text{if $\mu=\gamma$}, \\
0, & \text{otherwise},
\end{cases}$$
the multiplicity  (\ref{multiplicity}) is
$$
\begin{cases}
c^{\lambda}_{\nu\,  \mu}, & \text{if $\ell(\mu)\leq n$}, \\
0, & \text{otherwise}.
\end{cases}$$
Since $c^{\lambda}_{\nu \,\mu}=c^{\lambda}_{\mu \nu}$, we
have
$\textsf{s}^+_{\mu'}\left(z_{\{\lambda\}}\right)=z_{\{\lambda/\mu\}}$
when $\ell(\mu)\leq n$ and $0$ otherwise, by Lemma \ref{Jacobi for Ainfty}.
The proof for $\textsf{s}^-_{\mu'}(z_{\{\lambda\}})$ is similar.
\qed\vskip 3mm

\begin{cor} For $\mu\in\cP$ with $\ell(\mu)\leq n$, we have
\begin{equation*}
\text{\rm
$\textsf{s}^+_{\mu'}\left(z_{\{(0^n)\}}\right)$}=z_{\{\mu^\ast\}}, \ \ \
\text{\rm
$\textsf{s}^-_{\mu'}\left(z_{\{(0^n)\}}\right)$}=z_{\{\mu\}}.
\end{equation*}
\end{cor}
\vskip 3mm

For $n\geq 1$, we put
\begin{equation}
\begin{split}
&\textsf{h}^\pm_n=\textsf{s}^\pm_{(n)}=
\rho\left(t^\pm_{\{(n)\}}\right).
\end{split}
\end{equation}

\begin{prop}\label{hn action} Let $m\geq 0$ and $\lambda\in\Z_+^n$ be given.
\begin{itemize}
\item[(1)] If $m>n$, then {\rm $\textsf{h}^\pm_m$} acts as identically zero on $\R_n$.

\item[(2)] If $m\leq n$, then we have {\rm
\begin{equation*}
\begin{split}
\textsf{h}^+_m \left(z_{\{\lambda\}}\right)&=\sum_{\substack{\mu\in\Z^n_+ \\ (\lambda-(\mu_n^n))/(\mu-(\mu_n^n)) : \\
\text{a vertical strip of length $m$}}}z_{\{\mu\}},\\
\textsf{h}^-_m \left(z_{\{\lambda\}}\right)&=\sum_{\substack{\mu\in\Z^n_+ \\ (\mu-(\lambda_n^n))/(\lambda-(\lambda_n^n)) : \\
\text{a vertical strip of length $m$}}}z_{\{\mu\}}.
\end{split}
\end{equation*}
} In particular, we have {\rm $\textsf{h}^\pm_n
\left(z_{\{\lambda\}}\right)=z_{\{\lambda\mp(1^n)\}}$}.
\end{itemize}

\end{prop}
\pf (1) It follows directly from Theorem \ref{rho s lambda}.

(2) By Lemma \ref{Jacobi for Ainfty} and Theorem \ref{rho s lambda},
we have
$$\textsf{h}^-_m \left(z_{\{\lambda\}}\right)=z_{\{\lambda/-(0^{n-m},1^m)\}}=\sum_{\mu\in\Z_+^n}c^{\lambda}_{-(0^{n-m},1^m) \mu}z_{\{\mu\}}.$$
Since $c^{\lambda}_{-(0^{n-m},1^m)
\mu}=c^{\lambda+(1^n)}_{(1^{n-m},0^m)
\mu}=c^{\lambda+(1^n)-(\mu_n^n)}_{(1^{n-m},0^m) \mu-(\mu_n^n)}$, we
have
\begin{equation*}
\begin{split}
&c^{\lambda}_{-(0^{n-m},1^m) \mu} \\
&=
\begin{cases}
1, & \text{if $\lambda+(1^n)-(\mu_n^n)$ is a partition and  } \\
 & \ \ \ \text{$\left(\lambda+(1^n)-(\mu_n^n)\right)/\left(\mu-(\mu_n^n)\right)$ is a vertical strip of length $n-m$}, \\
0, & \text{otherwise}.
\end{cases}
\end{split}
\end{equation*}
Therefore, $c^{\lambda}_{-(0^{n-m},1^m) \mu}=1$ if and only if
$(\mu-(\lambda_n^n))/(\lambda-(\lambda_n^n))$ is a vertical strip of
length $m$. The proof for $\textsf{h}^+_m
\left(z_{\{\lambda\}}\right)$ is similar. \qed\vskip 3mm

As a corollary, we have the following.
\begin{cor}\label{reduction} As operators on $\R_n$ $(n\geq 1)$, we have for $0\leq i\leq n$ {\rm
\begin{equation*}
\textsf{h}^+_n\textsf{h}^-_i=\textsf{h}^+_{n-i}.
\end{equation*}
} Here, we assume {\rm $\textsf{h}^\pm_{0}={\rm id}_{\R_n}$}.
\end{cor}\vskip 3mm

\subsection{The action of $\K_0$ on $\K^{\rm l.w.}_{-n}$ and the $(\gl_\infty,\gl_n)$-duality}
For $n\geq 1$, let $R^\circ(GL_n(\mathbb{C}))$ be the character ring of finite
dimensional polynomial representations of $GL_n(\mathbb{C})$, which
is isomorphic to the ring of symmetric polynomials in $x_{[n]}$.
Also, it is free commutative algebra generated by $e_k(x_{[n]})$ for
$1\leq k\leq n$. Let $R(GL_n(\mathbb{C}))$ be the character ring of
finite dimensional representations of $GL_n(\mathbb{C})$. Then
$R(GL_n(\mathbb{C}))$ is the ring of symmetric Laurent polynomials
in $x_{[n]}$, and it is the localization of
$R^\circ(GL_n(\mathbb{C}))$ with respect to the multiplicative
subset $\{\,e_n(x_{[n]})^m\,|\,m\geq 1\,\}$.

\begin{thm}\label{K0 module structure}  For $n\geq 1$, let $\mathscr{M}_n=\K_0 \cdot [\B(-n\Lambda_0)]$ be the $\K_0$-submodule of $\K^{\rm l.w.}_{-n}$ generated by
$[\B(-n\Lambda_0)]$.
\begin{itemize}
\item[(1)] We have
    $$\mathscr{M}_n=\bigoplus_{\lambda\in\Z_+^n}\Z[\B(-\Lambda_\lambda)].$$
   In particular, $\K^{\rm l.w.}_{-n}$ is the completion of $\mathscr{M}_n$.

\item[(2)] There exists a $\Z$-algebra isomorphism
 $$R(GL_n(\mathbb{C})) \longrightarrow \K_0/{\rm ann}_{\K_0}(\mathscr{M}_n),$$ where
  $s_{\lambda}(x_{[n]})$ is mapped to $\overline{[\cB_{(\lambda+(\ell^n))',(\ell^n)'}]}$ for $\lambda\in\Z_+^n$ with $\ell=\max\{-\lambda_n,0\}$.

\item[(3)] The ideal ${\rm ann}_{\K_0}(\mathscr{M}_n)$ is generated by
$$[\cB_{(k),\emptyset}], [\cB_{\emptyset,(k)}]\ (k>n), \ \ [\cB_{(n),(i)}]-[\cB_{(n-i),\emptyset}]\ (0\leq i\leq n).
$$
\end{itemize}
\end{thm}
\pf (1) Since we may identify $\K_0$ with $\Z[t^+,t^-]$, and
$\mathscr{M}_n$ with $\Z[t^+,t^-]z_{\{(0^n)\}}$ in $\R_n$, it
is sufficient to show that $\mathscr{M}_n=\bigoplus_{\lambda\in\Z_+^n}\Z
z_{\{\lambda\}}$. By Theorem \ref{rho s lambda}, we
have $\mathscr{M}_n\subset\bigoplus_{\lambda\in\Z_+^n}\Z
z_{\{\lambda\}}$. Conversely, given $\lambda\in\Z_+^n$, put
$\mu=\lambda+(\ell^n)\in\cP$, where $\ell=\max\{-\lambda_n,0\}$. By
Theorem \ref{rho s lambda} and Proposition \ref{hn action}, we have
\begin{equation}\label{cyclic action}
z_{\{\lambda\}}=\left(\left(\textsf{h}^+_{n}\right)^\ell\circ
\textsf{s}^-_{\mu'}\right)\left(z_{\{(0^n)\}}\right)=\left(t^+_{\{(n)\}}\right)^\ell\cdot\left(
t^-_{\{\mu'\}}\cdot z_{\{(0^n)\}}\right),
\end{equation}
and hence $\bigoplus_{\lambda\in\Z_+^n}\Z
z_{\{\lambda\}}\subset\mathscr{M}_n$.

(2) Note that the map $g : \K_0/{\rm ann}_{\K_0}(\mathscr{M}_n) \rightarrow \mathscr{M}_n$ given by $g(\ov{a})=a\cdot z_{\{(0^n)\}}$ is a $\Z$-linear isomorphism. Define a ring homomorphism $f : R(GL_n(\mathbb{C})) \rightarrow
\K_0/{\rm ann}_{\K_0}(\mathscr{M}_n)$ by
$f(e_k(x_{[n]}))=\overline{t^-_{\{(k)\}}}$ ($1\leq k\leq n$) and
$f(e_n(x_{[n]})^{-1})=\overline{t^+_{\{(n)\}}}$. Since
$\overline{t^+_{\{(n)\}}t^-_{\{(n)\}}}=1$ in $\K_0/{\rm
ann}_{\K_0}(\mathscr{M}_n)$ by Corollary \ref{reduction}, $f$ is
well-defined. For $\lambda\in\Z_+^n$, we have
$$f(s_{\lambda}(x_{[n]}))
=f\left(e_n(x_{[n]})^{-\ell}s_{\lambda+(\ell^n)}(x_{[n]})\right)
=\overline{\left(t^+_{\{(n)\}}\right)^\ell\,
t^-_{\{(\lambda+(\ell^n))'\}}},$$ where $\ell=\max\{-\lambda_n,0\}$.
By (\ref{cyclic action}), we have
$g\circ f(s_{\lambda}(x_{[n]}))=z_{\{\lambda\}}$,
and hence obtain a $\Z$-linear isomorphism
\begin{equation*}
g\circ f : R(GL_n(\mathbb{C})) \longrightarrow \mathscr{M}_n.
\end{equation*}
This implies that $f$ is an isomorphism.

(3) Let $I_n=\langle\,t^\pm_{\{(m)\}}\,(m>n),\,
t^+_{\{(n)\}}t^-_{\{(i)\}}=t^+_{\{(n-i)\}}\,(0\leq i\leq
n)\,\rangle$, which is an ideal in $\Z[t^+,t^-]$. By Proposition
\ref{hn action} and Corollary \ref{reduction}, $I_n\subset {\rm
ann}_{\K_0}(\mathscr{M}_n)$. Hence we have an algebra homomorphism
$\iota : \K_0/I_n \longrightarrow \K_0/{\rm
ann}_{\K_0}(\mathscr{M}_n)$. Similarly, we have a surjective algebra
homomorphism $h : R(GL_n(\mathbb{C})) \rightarrow \K_0/I_n$ such
that $h(e_k(x_{[n]}))=\overline{t^-_{\{(k)\}}}$ ($1\leq k\leq n$)
and $h(e_n(x_{[n]})^{-1})=\overline{t^+_{\{(n)\}}}$. Since
$\iota\circ h=f$ and $f$ is an isomorphism by (2), it follows that $h$ is an isomorphism and hence so is $\iota$. This implies that $I_n = {\rm ann}_{\K_0}(\mathscr{M}_n)$ \qed\vskip
3mm

By Theorem \ref{K0 module structure}, $\mathscr{M}_n$ is naturally equipped with an $R(GL_n(\mathbb{C}))$-module structure, while $R(GL_n(\mathbb{C}))$ is itself an $R(GL_n(\mathbb{C}))$-module by left multiplication. Hence, we obtain the following.

\begin{cor}\label{duality via adjoints}
$\mathscr{M}_n$ is a free $R(GL_n(\mathbb{C}))$-module of rank 1, where
$s_{\lambda}(x_{[n]})$ corresponds to $[\B(-\Lambda_\lambda)]$ for $\lambda\in\Z_+^n$.
\end{cor}

\begin{rem}{\rm
Consider a graded $\Z$-module $R(GL)=\bigoplus_{n\geq 0}R(GL_n(\mathbb{C}))$, where $R(GL_0(\mathbb{C}))=\Z$. Then $R(GL)$ is a graded coalgebra with a graded comultiplication $\Delta$ given by $$\Delta(\chi)=\sum_{p+q=n}{\rm Res}^{GL_n(\mathbb{C})}_{GL_p(\mathbb{C})\times GL_q(\mathbb{C})}(\chi)$$ for $\chi \in R(GL_n(\mathbb{C}))$. Then by the $(\gl_\infty,\gl_n)$-duality (see Proposition \ref{duality-1}), we have as $\Z$-algebras
$$\bigoplus_{n\geq 0}\K^{\rm l.w.}_{-n} \simeq R(GL)^\circ,$$
where $R(GL)^\circ$ is the restricted dual of $R(GL)$ \cite{W99'}. Equivalently, the branching rule in $R(GL)$ corresponds to
the tensor product rule of integrable lowest weight $U_q(\gl_\infty)$-modules. On the other hand,
Theorem \ref{K0 module structure} explains this duality from a
different point of view. That is, we can recover $R(GL)$ via the
action of $\K_0$ on $\K^{\rm l.w.}$,
$$R(GL)=\Z\oplus\left( \bigoplus_{n\geq 1}\K_0/{\rm ann}_{\K_0}\left([\B(-n\Lambda_0)]\right)\right).$$
This explains that  the multiplication in $R(GL_n(\mathbb{C}))$,
which is not a graded multiplication in $R(GL)$, corresponds to the
composite of operators  on $\K^{\rm l.w.}_{-n}$ associated with
level zero extremal weight crystals. }
\end{rem}

\subsection{Hall-Littlewood vertex operators on $\K^{\rm l.w.}$}
Fix an indeterminate $q$. Let $\mathscr{B}^q$ be the associative $\Z[q]$-algebra generated by $\{\,B^q_n\,|\,n\in\Z\,\}$ subject to the relation;
\begin{equation}\label{defining rel for Bt}
B^q_mB^q_n=qB^q_nB^q_m+qB^q_{m+1}B^q_{n-1}-B^q_{n-1}B^q_{m+1} \ \ \ (m,n\in\Z).
\end{equation}
The above defining relations can be written in terms of generating functions as follows;
\begin{equation}
(u-qv){B}^q(u){B}^q(v)=(qu-v){B}^q(v){B}^q(u),
\end{equation}
where $u,v$ are formal commuting variables and
$B^q(u)=\sum_{k\in\Z}B^q_ku^k$. Note that $\mathscr{B}^q$ is
isomorphic to the ``half'' of the algebra generated by the
Hall-Littlewood vertex operators on $\Z[q]\otimes Sym$ introduced by
Jing \cite{J} (see also \cite{G,SZ}).

Our main claim in this section is that there is a natural embedding
of $\mathscr{B}^q$ into ${\rm End}(\K^{\rm l.w.})[[q]]$ and hence we can define a
Hall-Littlewood vertex operator for $\gl_\infty$. Put
\begin{equation}
\begin{split}
&\texttt{F}^\vee(u)=\sum_{k\in\Z} \rho\left([\B(-\Lambda_k)]\right)u^k,  \\
&\texttt{S}(v)=\sum_{n\geq 0 }\rho\left( [\cB_{(n)}]\right) v^n,\ \ \  \texttt{E}(v)=\sum_{n\geq 0 }\rho\left( [\cB_{(1^n)}]\right) v^n,
\end{split}
\end{equation}
where $\rho : \K^\vee \rightarrow {\rm End}(\K^{\rm l.w.})$ is given in (\ref{rho}).
We define
\begin{equation}
\begin{split}
\texttt{B}^q(u)&=\texttt{F}^\vee(u)\texttt{S}(-u^{-1})\texttt{E}(qu^{-1})=\sum_{k\in\Z}\texttt{B}^q_k u^k\in \left( {\rm End}(\K^\vee)[[q]] \right)[[u,u^{-1}]]. \\
\end{split}
\end{equation}
Equivalently, for $k\in\Z$
\begin{equation}\label{Bqk formula}
\begin{split}
\texttt{B}^q_k&=\sum_{i\geq 0}\sum_{j=0}^i(-1)^{i-j} q^j\rho\left(\left[\B(-\Lambda_{i+k})\otimes \cB_{(i-j)}\otimes \cB_{(1^j)}\right]\right) \\
&=\sum_{j\geq 0}\left(\sum_{i\geq 0}(-1)^{i} \rho\left(\left[\B(-\Lambda_{i+j+k})\otimes \cB_{(i)}\otimes \cB_{(1^j)}\right]\right)\right)q^j.
\end{split}
\end{equation}
Note that the coefficient of $q^j$ ($j\geq 0$) in $B^q_k$ gives a well-defined $\Z$-linear map  from $\K^{\rm l.w.}_{-n}$ to $\K^{\rm l.w.}_{-n-1}$ for $n\geq 0$ by Proposition \ref{hn action} (1), and hence a $\Z$-linear operator on $\K^{\rm l.w.}$.  

\begin{lem}\label{Pieri-general-2} We have {\rm
$$\texttt{S}(v)\texttt{F}^\vee(u) = \texttt{F}^\vee(u) \texttt{S}(v){(1+uv)}.$$}
\end{lem}
\pf  By Lemma \ref{Pieri-general}, we have for $i\in\Z$ and $k\geq 0$,
\begin{equation*}
\begin{split}
&\B(\Lambda_i)\otimes\cB_{(k)}^\vee\simeq
\left(\cB_{(k)}^\vee\otimes \B(\Lambda_{i})\right) \sqcup
\left(\cB_{(k-1)}^\vee\otimes \B(\Lambda_{i-1})\right).
\end{split}
\end{equation*} By taking its dual, we obtain the identity.
\qed\vskip 3mm

\begin{lem}\label{Bt relations} We have{\rm
$$(u-qv)\texttt{B}^q(u)\texttt{B}^q(v)=(qu-v)\texttt{B}^q(v)\texttt{B}^q(u).$$}
\end{lem}
\pf  It follows directly from the dual identity of (\ref{commuting relation-2}) and Lemma \ref{Pieri-general-2}.
\qed\vskip 3mm

By Lemma \ref{Bt relations}, $\texttt{B}^q_k$'s satisfy the defining relation
(\ref{defining rel for Bt}) for $\mathscr{B}^q$, and we have a $\Z[q]$-algebra homomorphism
\begin{equation}\label{t embedding}
\varpi : \mathscr{B}^q \longrightarrow  {\rm End}(\K^{\rm l.w.})[[q]]
\end{equation}
given by $\varpi(B^q_k)=\texttt{B}^q_k$ for $k\in\Z$.

We define  for $n\geq 1$ and $\alpha=(\alpha_1,\ldots,\alpha_n)\in\Z^n$
\begin{equation}\label{definition of Balphat}
{B}^q_\alpha = \prod_{1\leq i<j\leq n}(1-qR_{ij}){B}^q_{\alpha_1}\cdots {B}^q_{\alpha_n},
\end{equation}
where $R_{ij}$ is a raising operator, i.e.
$$R_{ij}\left({B}^q_{\alpha_1}\cdots
{B}^q_{\alpha_n}\right)={B}^q_{\alpha_1}\cdots{B}^q_{\alpha_i+1}\cdots{B}^q_{\alpha_j-1}\cdots
{B}^q_{\alpha_n}.$$ We should remark that when ${B}^q_k$'s are
replaced by the Hall-Littlewood vertex operators on $\Z[q]\otimes Sym$, ${B}^q_\lambda$
($\lambda\in\Z_+^n$) gives the operator introduced by Shimozono and
Zabrocki \cite{SZ}. By (\ref{defining rel for Bt}), we can check that
${B}^q_\alpha=-{B}^q_{(\alpha_1,\cdots,\alpha_{i+1}-1,\alpha_i+1,\cdots,\alpha_n)}$.
More generally, for each permutation $w$ in $S_n$, we have
\begin{equation}\label{shifted alternating}
{B}^q_\alpha=(-1)^{\ell(w)}{B}^q_{w(\alpha+\rho_n)-\rho_n},
\end{equation}
where $\ell(w)$ is the length of $w$ and
$\rho_n=(n-1,n-2,\ldots,0)$. Now for $\lambda\in\Z_+^n$, we put
\begin{equation}\label{Bqlambda}
\texttt{B}^q_\lambda=\varpi(B^q_\lambda).
\end{equation}
Then we have the following generating function for $\texttt{B}^q_\lambda$.
\begin{prop} For $n\geq 1$, we have {\rm
\begin{equation*}
\begin{split}
&\sum_{\lambda\in\Z_+^n}\texttt{B}^q_{\lambda}s_\lambda(x_{[n]})
=\texttt{F}^\vee(x_{[n]})
\texttt{S}(-x^{-1}_{[n]})\texttt{E}(qx^{-1}_{[n]}).
\end{split}
\end{equation*}}
\end{prop}
\pf By (\ref{commuting relation-2}) and Lemma \ref{Pieri-general-2}, we have
\begin{equation}\label{B product}
\begin{split}
&\texttt{B}^q(x_1)\cdots \texttt{B}^q(x_n)
\\
&=\prod_{1\leq i<j\leq n}\frac{1-x_i^{-1}x_j}{1-qx_i^{-1}x_j}\ \texttt{F}^\vee(x_{[n]})
\texttt{S}(-x^{-1}_{[n]})\texttt{E}(qx^{-1}_{[n]}).
\end{split}
\end{equation}
Since $$\prod_{1\leq i<j\leq n}(1-qx_i^{-1}x_j)\ \texttt{B}^q(x_1)\cdots \texttt{B}^q(x_n)
=\sum_{\lambda\in\Z_+^n}\texttt{B}^q_\lambda\left(\sum_{w\in S_n}(-1)^{\ell(w)}x_{[n]}^{w(\lambda+\rho_n)-\rho_n}\right)$$ by (\ref{shifted alternating}) and (\ref{Bqlambda}), we have
\begin{equation*}
\prod_{1\leq i<j\leq n}\frac{1-qx_i^{-1}x_j}{1-x_i^{-1}x_j}\
\texttt{B}^q(x_1)\cdots
\texttt{B}^q(x_n)=\sum_{\lambda\in\Z_+^n}\texttt{B}^q_{\lambda}s_\lambda(x_{[n]}).
\end{equation*} \qed\vskip 3mm

\begin{cor}\label{formula for Btlambda}
For $\lambda\in\Z_+^n$, we have {\rm
\begin{equation*}
\texttt{B}^q_\lambda=\sum_{\eta\in\Z_+^n}\sum_{\substack{\sigma,\mu,\nu\in\cP
\\  \ell(\sigma),\ell(\mu),\ell(\nu)\leq n }}
(-1)^{|\mu|}q^{|\nu|}c^{\lambda}_{\eta\,\sigma^*}c^{\sigma}_{\mu\,\nu}\rho\left(\left[\B(-\Lambda_\eta)\otimes
\cB_{\mu}\otimes\cB_{\nu'} \right]\right).
\end{equation*}}
\end{cor}
\pf Note that
\begin{equation*}
\texttt{S}(x_{[n]})=\sum_{\mu \in \cP}\rho\left( [\cB_{\mu}]\right) s_{\mu}(x_{[n]}),
\end{equation*}
which can be verified by modifying the arguments in Proposition
\ref{duality on En} using matrices with non-negative integral
entries (see for example, \cite{K07}). By Lemma \ref{Cauchy
identities}, we have
\begin{equation*}
\begin{split}
&\texttt{F}^\vee(x_{[n]})
\texttt{S}(-x^{-1}_{[n]})\texttt{E}(qx^{-1}_{[n]}) \\
&=\sum_{\eta\in\Z_+^n}\sum_{\substack{\mu,\nu\in\cP \\ \ell(\mu),\ell(\nu)\leq n }}(-1)^{|\mu|} q^{|\nu|}
\rho\left(\left[\B(-\Lambda_\eta)\otimes \cB_{\mu}\otimes\cB_{\nu'} \right]\right)s_{\eta}(x_{[n]})s_{\mu}(x_{[n]}^{-1})s_{\nu}(x_{[n]}^{-1}).
\end{split}
\end{equation*}
Note that $s_{\zeta}(x_{[n]}^{-1})=s_{\zeta^*}(x_{[n]})$  and
$c^{\sigma^*}_{\mu^*\,\nu^*}=c^{\sigma}_{\mu\,\nu}$ for
$\zeta,\sigma,\mu,\nu\in\Z_+^n$. Comparing the coefficient of
$s_\lambda(x_{[n]})$ on both sides, we obtain the identity.
\qed\vskip 3mm

\begin{prop}\label{Vertex operator}
The map $\varpi : \mathscr{B}^q\longrightarrow {\rm End}(\K^{\rm l.w.})[[q]]$ in
(\ref{t embedding}) is injective.
\end{prop}
\pf 
By (\ref{definition of Balphat}),  we have for
$\alpha=(\alpha_1,\ldots,\alpha_n)\in\Z^n$
\begin{equation*}
{B}^q_{\alpha_1}\cdots {B}^q_{\alpha_n} = \prod_{1\leq i<j\leq n}(1-q R_{ij})^{-1}{B}^q_\alpha,
\end{equation*}
where $R_{ij}$ acts on ${B}^q_\alpha$ by
$R_{ij}({B}^q_\alpha)={B}^q_{(\alpha_1,\cdots,\alpha_i+1,\cdots,\alpha_j-1,\cdots,\alpha_n)}$.
By the same arguments as in \cite[Ex.III.6.4]{Mac95}, we have
\begin{equation}\label{Kostka-2}
(-1)^{\ell(w)}{B}^q_{\alpha_1}\cdots {B}^q_{\alpha_n}=\sum_{\lambda\in\Z^n_+}K_{\lambda\,\mu}(q)B^q_{\lambda},
\end{equation}
where $w\in S_n$ and $\mu\in\Z_+^n$ are given by
$w(\alpha+\rho_n)-\rho_n=\mu$,  and $K_{\lambda\,\mu}(q)$ is the coefficient of $x_1^{\mu_1}\cdots x_n^{\mu_n}$ in ${\rm det}(x_i^{\lambda_i-i+j})/\prod_{1\leq i<j\leq n}(1-qx_i^{-1}x_j)$, or the 
Kostka-Foulkes polynomial associated with $\lambda,\mu\in\Z_+^n$.
Note that the sum on the right hand side is not necessarily finite.

Suppose that $b\in\mathscr{B}^q$ is given such that $\varpi(b)=0$.
By (\ref{Kostka-2}), we may write
\begin{equation*}
b=\sum_{n=1}^N\sum_{\lambda\in\Z_+^n}c_{\lambda}(q)B^q_\lambda
\end{equation*}
for some $N\geq 1$ and $c_\lambda(q)\in\Z[q]$, and
$$\varpi(b)=\sum_{n=1}^N\sum_{\lambda\in\Z_+^n}c_{\lambda}(q)\texttt{B}^q_\lambda.$$ 

We first assume that $c_{\mu}(0)\neq 0$ for some $\mu\in\Z_+^n$.
Then we have
\begin{equation*}
\varpi(b)\cdot 1
|_{q=0}=\sum_{n=1}^N\sum_{\lambda\in\Z_+^n}c_{\lambda}(0)\texttt{B}^0_\lambda
\cdot 1 =\sum_{n=1}^N\sum_{\lambda\in\Z_+^n}c_{\lambda}(0)z_\lambda,
\end{equation*}
since $[\B_{(k)}]\cdot 1=[\B_{(1^k)}]\cdot 1=0$ for $k> 0$. Since
$\{\,z_\lambda\,|\,\lambda\in\Z_+^n, n\geq 1\,\}\cup\{1\}$ is linearly independent, we
have $c_\lambda(0)=0$ for all $\lambda$, which is a
contradiction.  Now suppose that $c_\lambda (0)=0$ for all
$\lambda$.  If $b\neq 0$, then we have $b=q^mb'$ for some $m\geq 1$
and $b'\in \mathscr{B}^q$,  where $c_{\mu}(q)q^{-m}\in\Z[q]$ has a
non-zero constant term for some $\mu$. Since $ {\rm End}(\K^{\rm l.w.})[[q]]$
is free over $\Z[q]$, we have $\varpi(b')=0$. By the same argument,
we conclude that $b'=0$, and hence $b=0$, which is also a
contradiction. Therefore, $\varpi$ is injective. \qed\vskip
3mm

\begin{prop}\label{Hall Littlewood} For $\mu\in\Z_+^n$, we have {\rm
\begin{equation*}
\texttt{B}^q_{\mu_1}\cdots \texttt{B}^q_{\mu_n}\cdot 1=
\frac{1}{\prod_{1\leq i<j\leq
n}(1-qR_{ij})}z_{\{\mu\}}=\sum_{\lambda\in\Z^n_+}K_{\lambda\mu}(q)z_{\{\lambda\}},
\end{equation*}}
where $R_{ij}$ acts on $z_{\{\mu\}}$ by
$R_{ij}\left(z_{\{\mu\}}\right)=z_{\{R_{ij}(\mu)\}}$.
\end{prop}
\pf By (\ref{B product}), we have
\begin{equation*}
\begin{split}
&\texttt{B}^q(x_1)\cdots \texttt{B}^q(x_n)
\\
&=\prod_{1\leq i<j\leq n}\frac{1-x_i^{-1}x_j}{1-qx_i^{-1}x_j}\ \texttt{F}^\vee(x_{[n]})
\texttt{S}(-x^{-1}_{[n]})\texttt{E}(qx^{-1}_{[n]}).
\end{split}
\end{equation*}
If we apply both sides to $1=[\B(0)]\in \K^{\rm l.w.}$, then
\begin{equation*}
\begin{split}
&\texttt{B}^q(x_1)\cdots \texttt{B}^q(x_n)\cdot 1
=\prod_{1\leq i<j\leq n}\frac{1-x_i^{-1}x_j}{1-qx_i^{-1}x_j}\texttt{F}^\vee(x_{[n]})\cdot 1,
\end{split}
\end{equation*}
since $[\cB_{(k)}]\cdot 1=[\cB_{(1^k)}]\cdot 1=  0$ for $k>0$.
Given $\mu\in\Z_+^n$, equating the coefficients of $x_1^{\mu_1}\cdots x_n^{\mu_n}$ on both sides, we obtain the first identity since $z_{\{\mu\}}=\prod_{1\leq i<j\leq n}(1-R_{ij})z_\mu$. The second identity follows from the same arguments as in (\ref{Kostka-2}).
\qed

\begin{rem}{\rm
If $q=0$, then $\mathscr{B}^0$ is the algebra of Bernstein
operators \cite{Mac95} and $\texttt{B}^0_{\lambda_1}\cdots
\texttt{B}^0_{\lambda_n}\cdot
1=z_{\{\lambda\}}=[\B(-\Lambda_\lambda)]$ ($\lambda\in\Z_+^n$),
which is a Rodrigues type formula. Also, when $q=1$, we have
$\texttt{B}^1_{\lambda_1}\cdots \texttt{B}^1_{\lambda_n}\cdot
1=z_{\lambda}=[\B(-\Lambda_{\lambda_1})\otimes\cdots\otimes
\B(-\Lambda_{\lambda_n})]$. Hence $\texttt{B}^q_{\lambda_1}\cdots
\texttt{B}^q_{\lambda_n}\cdot 1$ interpolates $z_{\{\lambda\}}$ and
$z_\lambda$, and $\texttt{B}^q_k$ may be viewed as a
Hall-Littlewood vertex operator for $\gl_\infty$. }
\end{rem}\vskip 3mm

One may define another Hall-Littlewood vertex operator for $\gl_{\infty}$ using the involution  $\omega$  on $\cD$  (\ref{omega}). 
Put
\begin{equation}
\overline{\texttt{B}}^q(u)=\sum_{k\in\Z}\overline{\texttt{B}}^q_k u^k=\texttt{F}^\vee(u^{-1})\texttt{S}^\vee(-u^{-1})\texttt{E}^\vee(qu^{-1}), 
\end{equation}
where
\begin{equation*}
\begin{split}
&\texttt{S}^\vee(v)=\sum_{n\geq 0 }\rho\left( [\cB^\vee_{(n)}]\right) v^n,\ \ \  \texttt{E}^\vee(v)=\sum_{n\geq 0 }\rho\left( [\cB^\vee_{(1^n)}]\right) v^n.
\end{split}
\end{equation*}
Note that for $k\in\Z$ 
\begin{equation}
\begin{split}
\ov{\texttt{B}}^q_k&=\sum_{j\geq 0}\left(\sum_{i\geq 0}(-1)^{i} \rho\left(\left[\B(-\Lambda_{-i-j-k})\otimes \cB^\vee_{(i)}\otimes \cB^\vee_{(1^j)}\right]\right)\right)q^j \\
&=\sum_{j\geq 0}\left(\sum_{i\geq 0}(-1)^{i} \rho\left(\omega\left[\B(-\Lambda_{i+j+k})\otimes \cB_{(i)}\otimes \cB_{(1^j)}\right]\right)\right)q^j.
\end{split}
\end{equation}
By (\ref{Bqk formula}), we also have a $\Z[q]$-algebra homomorphism $\ov{\varpi} : \mathscr{B}^q \rightarrow {\rm End}(\K^{\rm l.w.})[[q]]$ given by $\ov{\varpi}(B^q_k)=\ov{\texttt{B}}^q_k$ for $k\in \Z$.
 
\begin{prop} For $m,n\in\Z$, we have{\rm
\begin{equation*}
\overline{\texttt{B}}^q_m \texttt{B}^q_n = \texttt{B}^q_n \overline{\texttt{B}}^q_m.
\end{equation*}}
\end{prop}
\pf It is straightforward to check that
\begin{equation*}
\begin{split}
&{\texttt{B}}^q(u)\overline{\texttt{B}}^q(v) \\ &=\texttt{F}^\vee(u)\texttt{S}(-u^{-1})\texttt{E}(qu^{-1})\texttt{F}^\vee(v^{-1})\texttt{S}^\vee(-v^{-1})\texttt{E}^\vee(qv^{-1}) \\
&=\texttt{F}^\vee(u)\texttt{F}^\vee(v^{-1})\texttt{S}(-u^{-1})\texttt{E}(qu^{-1})\texttt{S}^\vee(-v^{-1})\texttt{E}^\vee(qv^{-1})\frac{(1-u^{-1}v^{-1})}{(1-qu^{-1}v^{-1})} \\
&=\overline{\texttt{B}}^q(v){\texttt{B}}^q(u),
\end{split}
\end{equation*}
which implies that $\overline{\texttt{B}}^q_m \texttt{B}^q_n = \texttt{B}^q_n \overline{\texttt{B}}^q_m$ for $m,n\in\Z$.
\qed\vskip 3mm
 
\begin{rem}{\rm
For $\lambda\in\Z_+^n$, we have $\ov{\texttt{B}}^0_{\lambda_1}\cdots
\ov{\texttt{B}}^0_{\lambda_n}\cdot
1=z_{\{\lambda^*\}}=[\B(-\Lambda_{\lambda^*})]$.
}
\end{rem}

{\small

\end{document}